\documentclass[a4paper,11pt]{article}
\input epsf
\usepackage{latexsym}
\usepackage{amssymb}
\newtheorem{thm}{Th\'eor\`eme}[section]
\newtheorem{lem}[thm]{Lemme}

\newtheorem{cor}[thm]{Corollaire}
\newtheorem{prop}[thm]{Proposition}

\newtheorem{rem}[thm]{Remarque}

\input epsf

\title{Sur le groupe d'interpolation}
\author{Roland Bacher}
\date{}

\begin{document}
\maketitle


{\it R\'esum\'e : Nous \'etudions le groupe d'interpolation
dont les \'el\'ements
sont certaines paires de s\'eries formelles. Ce groupe poss\`ede
une repr\'esentation lin\'eaire fid\`ele dans les matrices
triangulaires inf\'erieures infinies. On peut donc le munir
d'une structure de groupe de Lie naturelle. Une fonction 
reli\'ee \`a l'exponentielle matricielle entre son alg\`ebre de Lie 
et sa repr\'esentation lin\'eaire \'etend alors 
l'exponentielle usuelle \`a deux arguments (qui sont des s\'erie
formelles) et cette extension poss\`ede des propri\'et\'es 
int\'eressantes. Nous terminons par une application \`a la
combinatoire \'enum\'erative et la description d'une alg\`ebre qui
g\'en\'eralise le groupe d'interpolation.}

{\it Abstract: 
We study the interpolation group whose elements are suitable pairs of
formal power series. This group has a faithful
representation into infinite lower triangular matrices
and carries thus a natural structure as a Lie group. 
The matrix exponential between its Lie algebra and 
its matrix representation gives rise to a function
with interesting properties  
extending the usual exponential function to two variables
(which are formal power series) We finish with an application to
enumerative combinatorics and the description of an algebra 
which generalizes the interpolation group.}

\section{Introduction}

Notons $\mathfrak m=x\mathbb C[[x]]$ l'id\'eal maximal des s\'eries 
formelles sans terme constant et consid\'erons le groupe
multiplicatif $\mathcal U=\mathbb C[[x]]\setminus \mathfrak m$ 
des germes inversibles de fonctions au voisinage de $0\in\mathbb C$.
Introduisons \'egalement le groupe non-commutatif
$\mathcal D=\mathfrak m\setminus
\mathfrak m^2$ des ``diff\'eomorphismes formelles'' (avec loi de
groupe donn\'ee par la composition des s\'eries).
L'action \'evidente de $\mathcal D$ sur $\mathcal U$ par
``changements de cartes formelles'' permet de
consid\'erer le produit semi-direct 
$\mathcal I=\mathcal U\rtimes \mathcal D$ qu'on appellera le {\it
groupe d'interpolation}. Notre but est de d\'ecrire
quelques propri\'et\'es de ce groupe. Ce groupe (ou plut\^ot 
sa repr\'esentation matricielle) a \'egalement \'et\'e introduit
sous le nom de ``groupe de Riordan'' par L. Shapiro dans le but
d'interpr\'eter certains triangles de nombres combinatoires, voir par 
exemple \cite{SGWW}.

Le reste de ce papier 
est organis\'e comme suit :

Dans le chapitre \ref{sectgroupeI} nous consid\'erons 
quelques propri\'et\'es \'el\'ementaires de $\mathcal I$.

Le sous-groupe $\mathcal{SI}$ obtenu \`a partir des sous-groupes
$\mathcal{SU}=1+\mathfrak m\subset \mathcal U$ et $\mathcal{SD}=
x+\mathfrak m^2\subset \mathcal D$ intervient dans le chapitre 
\ref{sectinterpol} pour d\'ecrire deux d\'eformations ``naturelles''
(et holomorphes sur les sous-groupes form\'es de s\'eries 
holomorphes au voisinage de $0$) du groupe ab\'elien
$\mathcal{SU}$ vers le groupe non-ab\'elien $\mathcal{SD}$.

Le groupe $\mathcal I$ admet une repr\'esentation matricielle
fid\`ele $\rho$ dans les matrices triangulaires inf\'erieures
infinies. Cette repr\'esentation, d\'ecrite dans le
chapitre \ref{sectreprmatr}, permet de consid\'erer $\mathcal I$ 
comme un groupe de Lie de dimension infinie. Le chapitre 
\ref{sectalgLie} d\'ecrit l'alg\`ebre de Lie $\mathfrak i$
de $\rho(\mathcal I)\sim \mathcal I$.

Le chapitre \ref{sectexp} introduit une fonction $\hbox{Exp}:
\mathbb C[[x]]\times \mathbb C[[x]]
\sim\mathfrak i\longrightarrow \mathcal U$ 
obtenue en projettant l'application exponentielle usuelle
$\mathfrak i\longrightarrow \mathcal I=\mathcal U\rtimes \mathcal D$
sur le premier facteur. On peut consid\'erer
l'application $\hbox{Exp}$ comme une extension ou une d\'eformation
de l'exponentielle usuelle.

Les chapitres \ref{sectequadiff} et \ref{sectdevserie} d\'ecrivent 
des \'equations diff\'erentielles et un d\'eveloppement en s\'erie
pour $\hbox{Exp}(\alpha;\beta)\in\mathcal U$.

Dans le chapitre \ref{sectfctrecipr} nous \'etudions bri\`evement
la fonction r\'eciproque
de $\alpha \longmapsto \hbox{Exp}(\alpha;\beta)$ (pour $\beta
\in\mathbb C[[x]]$ fix\'e) et de $\beta\longmapsto
\hbox{Exp}(\beta;\beta)$.

Quelques propri\'et\'es analytiques de $\hbox{Exp}(\alpha;\beta)$
sont \'etudi\'ees dans le chapitre \ref{sectconvergnce}.

Le chapitre \ref{sectintcomb} donne une interpr\'etation
combinatoire pour la s\'erie de $\hbox{Exp}(s(x\beta)',s\beta)$ 
et utilise les r\'esultats du chapitre \ref{sectequadiff}
pour donner une nouvelle preuve d'un r\'esultat classique 
de combinatoire \'enum\'erative.

Le chapitre \ref{sectriordan} est un bref r\'esum\'e des travaux
effectu\'es par diff\'erents auteurs autour des ``matrices
de Riordan'' (les \'el\'ements de $\rho(\mathcal I)$). 

Finalement, le chapitre \ref{sectmatrpolyn}
d\'ecrit une alg\`ebre contenant un groupe matricielle
qui g\'en\'eralise le groupe $\rho(\mathcal{SI})$.

\begin{rem} \`A cause de l'ambiguit\'e de la notation
$f(gh)$, nous utiliserons toujours
$f\circ(gh)$ pour la composition de fonctions (ici $f$ compos\'e
avec le produit de $g$ et de $h$). Pour cette m\^eme raison,
nous remplacerons souvent des parenth\`eses par des accolades.
Ainsi $f\lbrace gh+1\}$ d\'esignera le produit de la fonction $f$ par
la
fonction obtenue en rajoutant $1$ au produit de $g$ et $h$.
\end{rem}  

\section{Le groupe d'interpolation $\mathcal I$}\label{sectgroupeI}

L'anneau commutatif ${\mathbb C}[[x]]$ des s\'eries formelles en $x$
est un anneau local d'id\'eal maximal ${\mathfrak m}=x{\mathbb
  C}[[x]]$ les s\'eries formelles sans terme constant.
D\'esignons par ${\mathcal U}={\mathbb C}[[x]]\setminus{\mathfrak m}=
{\mathbb C}^*+{\mathfrak m}$ le groupe
des unit\'es et par $\mathcal{SU}=1+{\mathfrak m}
\subset {\mathcal U}$
le sous-groupe correspondant aux s\'eries formelles avec coefficient
constant $1$. Le groupe $\mathcal{SU}$ est \'egalement
le noyau de l'homomorphisme $A\longmapsto A(0)$
qui associe \`a une s\'erie formelle $A=\sum_{n=0}^\infty A_nx^n
\in{\mathcal U}$ son coefficient constant $A(0)=A_0\in{\mathbb C}^*$. 
 
Notons de m\^eme ${\mathcal D}={\mathbb C}^* x+{\mathfrak m}^2$ 
l'ensemble des s\'eries formelles admettant une s\'erie r\'eciproque
(inverse pour la composition des s\'eries). On peut penser \`a
un \'el\'ement de ${\mathcal D}$ comme le germe d'un diff\'eomorphisme
formel en $0\in {\mathbb C}$. L'ensemble $\mathcal D$ est muni
de la structure de groupe non-commutatif $\alpha\beta=
\beta\circ \alpha$ donn\'e par la composition 
des s\'eries $\alpha,\beta\in{\mathcal D}$.
Nous \'ecrivons $\mathcal{SD}=x+x{\mathfrak m}=x+\mathfrak m^2$ 
pour le
sous-groupe des diff\'eomorphismes formels tangents \`a 
l'identit\'e. On peut \'egalement d\'efinir le sous-groupe 
$\mathcal{SD}\subset {\mathcal D}$ comme le noyau de 
l'homomorphisme $\alpha\longmapsto \alpha'(0)\in{\mathbb C}$
qui associe \`a $\alpha=\sum_{n=1}^\infty \alpha_nx^n\in{\mathcal D}$ 
sa d\'eriv\'ee $\alpha_1=\alpha'(0)\in{\mathbb C}^*$ \`a l'origine.

Le groupe ${\mathcal D}$
agit par automorphismes sur ${\mathcal U}$ en consid\'erant 
l'application $A\in{\mathcal U}\longmapsto 
A\circ \alpha\in {\mathcal U}$. Cette action se restreint en
une action de $\mathcal{D}$ (ou de son sous-groupe $\mathcal{SD}$)
sur $\mathcal{SU}$.

{\bf D\'efinition} Le {\it groupe d'interpolation 
$\mathcal  I$ } est
le produit semi-direct $\mathcal I={\mathcal U}\rtimes {\mathcal D}$
de ${\mathcal U}$ avec ${\mathcal D}$. Le groupe
d'interpolation sp\'ecial $\mathcal{SI}\subset{\mathcal I}$ 
est le produit semi-direct $\mathcal{SU}\rtimes \mathcal{SD}$.

Un \'el\'ement de $\mathcal I$ sera not\'e
$(A,\alpha),(B,\beta),(C,\gamma),\dots$ avec $A,B,C,\dots
\in{\mathcal U}$ et $\alpha,\beta,\gamma,\dots\in{\mathcal D}$.
Le produit dans $\mathcal I$ est donn\'e par
$$(A,\alpha)(B,\beta)=(A\lbrace B\circ \alpha\rbrace,\beta\circ \alpha)$$
pour $(A,\alpha),(B,\beta)\in I$.
L'inverse $(A,\alpha)^{-1}$ d'un \'el\'ement s'obtient par la formule
$$(A,\alpha)^{-1}=\left(\frac{1}{A\circ\alpha^{\langle -1\rangle}},
\alpha^{\langle -1\rangle}\right)$$
o\`u la s\'erie r\'eciproque $\alpha^{\langle
  -1\rangle}$ de $\alpha\in{\mathcal D}$ est d\'efinie par 
les identit\'es $\alpha\circ \alpha^{\langle -1\rangle}=
 \alpha^{\langle -1\rangle}\circ \alpha=x$. La structure de produit
semi-direct sur $\mathcal I$ \'equivaut \`a l'existence 
d'une suite exacte scind\'ee
$$0\longrightarrow {\mathcal U}\longrightarrow {\mathcal I}
={\mathcal U}\rtimes {\mathcal D}\longrightarrow {\mathcal D} 
\longrightarrow 1\ .$$
Pour l'injection ${\mathcal U}\longrightarrow {\mathcal I}$ 
nous choisirons dor\'enavant toujours
$A\in{\mathcal U}\longmapsto (A,x)\in \mathcal I$ ce qui nous
permettra  d'identifier $\mathcal U$ avec le sous-groupe
$(\mathcal U,x)\subset \mathcal I$. 
La surjection
${\mathcal I}\longrightarrow {\mathcal D}$ poss\`ede (par exemple)
la section
$\alpha\in{\mathcal D}\longmapsto (1,\alpha)\in {\mathcal I}$.
Nous noterons $(1,\mathcal D)\subset \mathcal I$ le sous-groupe
image de cette section.

Le groupe d'interpolation sp\'ecial $\mathcal{SI}=\mathcal{SU}
\rtimes\mathcal{SD}$ s'obtient comme le noyau de l'homomorphisme
${\mathcal I}\longrightarrow {\mathbb C}^*\times {\mathbb C}^*$ d\'efini
par l'\'evaluation
$(A,\alpha)\longmapsto (A(0),\alpha'(0))$. La suite exacte pour
$\mathcal I$ donne par restriction une suite exacte pour $\mathcal{SI}$.

\begin{prop} \label{propvarphi}
Pour $\kappa,\lambda,\mu\in{\mathbb C}$ trois nombres complexes,
l'application
$\varphi_{\kappa,\lambda,\mu}: 
\mathcal{SI}\longrightarrow \mathcal{SI}$ d\'efinie par
$$\varphi_{\kappa,\lambda,\mu}(A,\alpha)=\left(A^\kappa
\left\lbrace\frac{\alpha}{x}\right\rbrace^\lambda\left\lbrace
\alpha'\right\rbrace^\mu,\alpha\right)$$
est un endomorphisme de groupes. 
Pour $\kappa\in\mathbb C^*$, c'est un isomorphisme d'inverse
$\varphi_{1/\kappa,-\lambda,-\mu}$. 
\end{prop}

\begin{rem} \label{remrvtmntuniv} 
L'homomorphisme $\varphi_{\kappa,\lambda,\mu}$ est d\'efini
en utilisant la d\'etermination principale
$\hbox{log}\circ A,
\hbox{log}\circ\left(\frac{\alpha}{x}\right),
\hbox{log}\circ \alpha'
\subset\mathfrak m$ dans le calcul des puissances de
$A,\left(\frac{\alpha}{x}\right),\alpha'
\in 1+\mathfrak m$.

\`A cause de la monodromie du logarithme complexe,
l'homo\-morphisme $\varphi_{\kappa,\lambda,\mu}$ de la proposition 
\ref{propvarphi} ne provient pas d'un homomorphisme de $\mathcal I$,
sauf si $\kappa,\lambda,\mu\in\mathbb Z$.
On peut cependant l'\'etendre au produit semi-direct
${\mathcal U}_{>0}\rtimes {\mathcal D}_{>0}$
o\`u $\mathcal U_{>0}={\mathbb R}_{>0}+{\mathfrak m}\subset
{\mathcal U}$ et 
$\mathcal D_{>0}=x({\mathbb R}_{>0}+{\mathfrak m})\subset
{\mathcal D}$
ou encore
relever $\varphi_{\kappa,\lambda,\mu}$ sur le ``rev\^etement universel''
$\tilde{\mathcal I}$ de $\mathcal I$
obtenu par un rel\`evement dans $\mathbb R$ des arguments 
$\hbox{arg}(A(0)), \hbox{arg}(\alpha'(0))\in\mathbb R/(2\pi\mathbb Z)$.

Indiquons aussi l'existence des isomorphismes 
$$(A,\alpha)\longmapsto (A\lbrace A(0)\rbrace^a\lbrace
\alpha'(0)\rbrace^b,\alpha)$$
pour $a,b,\in\mathbb Z$ de $\mathcal I$ (sur le
groupe $\tilde\mathcal I$, ces homomorphismes sont d\'efinis
pour tout $a,b\in\mathbb C$). Leur nature est assez trivial
car le sous-groupe $(\mathbb C^*,x)\subset \mathcal I$ est le centre
de $\mathcal I$.
\end{rem}

\begin{rem} \label{remautomtriv} L'automorphisme int\'erieur
$$(A,\alpha)\longmapsto (1,Kx)(A,\alpha)\left(1,\frac{x}{K}\right)=
\left(A\circ(Kx),\frac{\alpha\circ (Kx)}{K}\right)$$
d\'efini pour tout $K\in\mathbb C^*$
correspond essentiellement \`a un changement de variable lin\'eaire
de $\mathbb C$.
\end{rem}

{\bf Preuve de la proposition \ref{propvarphi}} On a
$\varphi_{\kappa,\lambda,\mu}(1,x)=(1,x)$. Le calcul
$$\begin{array}{lcl}
\varphi_{\kappa,\lambda,\mu}(A,\alpha)\varphi_{\kappa,\lambda,\mu}(B,\beta)&=&
\left(A^\kappa\left\lbrace\frac{\alpha}{x}\right\rbrace^\lambda
\left\lbrace\alpha'\right\rbrace^\mu,\alpha\right)
\left(B^\kappa\left\lbrace\frac{\beta}{x}\right\rbrace^\lambda
\left\lbrace\beta'\right\rbrace^\mu,\beta\right)\\
&=&
\left(A^\kappa\left\lbrace\frac{\alpha}{x}\right\rbrace^\lambda
\left\lbrace\alpha'\right\rbrace^\mu
\left(B^\kappa\left\lbrace\frac{\beta}{x}\right\rbrace^\lambda
\left\lbrace\beta'\right\rbrace^\mu\right)\circ\alpha,
\beta\circ\alpha\right)\\
&=&
\left(A^\kappa\lbrace B\circ\alpha\rbrace^\kappa\left\lbrace
\frac{\beta\circ\alpha}
{x}\right\rbrace^\lambda
\left\lbrace(\beta\circ\alpha)'\right\rbrace^\mu,\beta\circ\alpha\right)
\\
&=&\varphi_{\kappa,\lambda,\mu}(
A\lbrace B\circ\alpha\rbrace,\beta\circ\alpha)\end{array}$$
termine la preuve.\hfill$\Box$

\begin{rem} On peut g\'en\'eraliser le groupe $\mathcal I$ (et son 
sous-groupe $\mathcal{SI}$) en rempla\c cant le groupe des
unit\'es ${\mathcal U}\subset {\mathbb C}[[x]]$ par le groupe
des unit\'es dans les germes de fonctions au voisinage d'un point
$P\in X$ o\`u $X$ est un espace topologique et en consid\'erant 
\`a la place de $\mathcal D$ un groupe form\'e de germes
d'hom\'eomorphismes fixant $P$.

Une autre g\'en\'eralisation de ${\mathcal I}$ consiste \`a remplacer
le corps ${\mathbb C}$ par un autre corps de base. Une grande partie
de notre papier s'adapte facilement au cas d'un corps 
commutatif quelconque.

Le groupe $\mathcal I$ poss\`ede des sous-groupes 
int\'eressants : On peut par exemple se restreindre aux s\'eries
ayant un rayon de convergence $>0$ et travailler uniquement avec
des fonctions $A\in{\mathcal U},\alpha\in\mathcal D$
holomorphes au voisinage de $0$. On peut \'egalement
se restreindre aux fonctions $A\in{\mathcal U},\alpha\in\mathcal D$
qui sont alg\'ebriques. (Pour cela, il faut montrer l'alg\'ebricit\'e
de la composition de deux fonctions alg\'ebriques et de la r\'eciproque
d'une fonction alg\'ebrique : Pour $\alpha,\beta$ alg\'ebriques dans
des ouverts convenables de ${\mathbb C}$, la composition 
$\gamma=\beta\circ\alpha$ v\'erifie l'\'equation 
$P(\alpha(x),\gamma(x))$ pour $P(x,\beta(x))=0$ une \'equation 
polynomiale d\'efinissant $\beta$. Le corps de fonctions
engendr\'e par $x$ et $\gamma$ est donc bien une extension finie
du corps des fonctions rationelles ${\mathbb C}(x)$ sur ${\mathbb
  C}$ sur la sph\`ere de Riemann. Similairement, on a
$P(\alpha^{\langle-1\rangle}(y),y)=0$ pour $P(x,\alpha(x))$
une \'equation polynomiale d\'efinissant la fonction alg\'ebrique 
$\alpha\in \mathcal D$.) 
\end{rem}

\begin{rem} La notation ${\mathcal U}$ choisie pour le groupe 
des unit\'es dans $\mathbb C[[x]]$ ne doit pas \^etre 
confondue avec la notation $\mathbf U({\mathcal H})$
utilis\'ee habituellement pour le groupe unitaire d'un espace de Hilbert.
\end{rem}

\section{Interpolations continues entre inversion et r\'e\-version}
\label{sectinterpol}

Pour $\tau\in\mathbb C$, consid\'erons les sous-ensembles
$$\mathcal{SG}(\tau)=\{(A,xA^\tau)\vert A\in\mathcal {SU}\}
\subset \mathcal{SI}$$
et
$$\mathcal{SG'}(\tau)=\{(A,\int_0 A^\tau)\vert A\in\mathcal{SU}\}
\subset \mathcal{SI}$$
o\`u l'on utilise la d\'etermination principale
$\hbox{log}(1+\mathfrak m)\subset \mathfrak m$ du logarithme
pour le calcul de $A^\tau=e^{\tau\hbox{log}(A)}$
et o\`u $\int_0A^\tau=x+\sum_{n=1}^\infty
A_{n,\tau}\frac{x^{n+1}}{n+1}$ 
est la primitive dans l'id\'eal maximal $\mathfrak
m=x\mathbb C[[x]]$ de la s\'erie $A^\tau=1+\sum_{n=1}^\infty
A_{n,\tau}x^n\in \mathcal{SU}$.

\begin{prop} \label{propIlassgr} $\mathcal{SG}(\tau)$ et
$\mathcal{SG'}(\tau)$ sont des sous-groupes de $\mathcal{SI}$, 
isomorphes \`a $\mathcal{SU}$ pour $\tau=0$
et isomorphes \`a $\mathcal{SD}$ sinon. 
\end{prop}

{\bf Preuve} La preuve est \'evidente pour $\tau=0$.
Pour $\tau\not=0$, elle r\'esulte des identit\'es
$$\mathcal{SG}(\tau)=\varphi_{0,\frac{1}{\tau},0}(1,\mathcal{SD})
=\{(\left\lbrace\frac{\alpha}{x}\right\rbrace^{1/\tau},\alpha),
\alpha\in\mathcal{SD}\}\subset \mathcal{SI}$$
et
$$\mathcal{SG'}(\tau)=
\varphi_{0,0,\frac{1}{\tau}}(1,\mathcal{SD})=\{(\left\lbrace\alpha'
\right\rbrace^{1/\tau},\alpha),\alpha\in\mathcal{SD}\}
\subset \mathcal{SI}$$
et de l'observation que $\varphi_{0,\lambda,\mu}(1,\mathcal{SD})$
est isomorphe \`a $\mathcal{SD}$ pour tout $\lambda,\mu
\in\mathbb C$.\hfill$\Box$

\begin{cor} \label{corinterpol} (i)
L'application 
$$\tau\longmapsto
\frac{1}{A\circ\left(xA^\tau\right)^{\langle -1\rangle}}$$
permet d'interpoler entre la s\'erie inverse 
$\frac{1}{A}$  (correspondante \`a $\tau=0$) de $A\in\mathcal{SU}$
et la s\'erie r\'eciproque 
$$(xA)^{\langle -1\rangle}=
\frac{x}{A\circ\left(xA\right)^{\langle -1\rangle}}$$
(correspondante \`a $\tau=1$) de
$xA\in\mathcal{SD}$.

\ \ (ii) L'application 
$$\tau\longmapsto
\frac{1}{A\circ\left(\int_0A^\tau\right)^{\langle -1\rangle}}\ ,$$
permet d'interpoler entre la s\'erie inverse 
$\frac{1}{A}$  (correspondante \`a $\tau=0$) de $A\in\mathcal{SU}$
et la s\'erie r\'eciproque 
$$(\int_0A)^{\langle -1\rangle}=
\int_0\frac{1}{A\circ\left(\int_0 A\right)^{\langle -1\rangle}}
$$ 
(correspondante \`a $\tau=1$) de $\int_0A\in\mathcal{SD}$. 
\end{cor}

Ce corollaire est \`a l'origine de la terminologie 
``groupe d'interpolation'' pour
$\mathcal{I}={\mathcal U}\rtimes {\mathcal D}$.

\begin{rem} Les deux interpolations du corollaire \ref{corinterpol}
se font de fa\c con
holomorphe (par rapport \`a $\tau$ et $x$) 
si $A\in \mathcal{SU}$ est holomorphe au voisinage de $0$.\end{rem}

{\bf Preuve du corollaire \ref{corinterpol}} Cela r\'esulte des
formules
$$\left(A,xA^\tau\right)^{-1}=
\left(
\frac{1}{A\circ \left(xA^\tau\right)^{\langle  -1\rangle}},
\left(xA^\tau\right)^{\langle -1\rangle}\right)\in\mathcal{SG}(\tau)\ ,$$
$$\left(A,\int_0A^\tau\right)^{-1}=
\left(
\frac{1}{A\circ \left(\int_0A^\tau\right)^{\langle  -1\rangle}},
\left(\int_0A^\tau\right)^{\langle -1\rangle}\right)\in
\mathcal{SG'}(\tau)$$
pour les \'el\'ements inverses de $(A,xA^\tau)\in\mathcal{SG}(\tau)$
et de
$(A,\int_0A^\tau)\in\mathcal{SG'}(\tau)$.
\hfill$\Box$

\begin{rem} Le calcul de
  $(C_\tau(s),\gamma_\tau(s))=(A,xA^\tau)^s\in\mathcal{SG}(\tau)\subset\mathcal{SI}$ (respectivement de 
$(\tilde C_\tau(s),\tilde \gamma_\tau(s))=(A,\int_0A^\tau)^s
\in\mathcal{SG'}(\tau)\subset\mathcal{SI}$)
 permet d'interpoler entre
$A^s=C_0(s)\in\mathcal{SU}$ et $(xA)^{\langle s\rangle}=xC_1(s)=\gamma_1(s)
\in \mathcal{SD}$ (respectivement $(\int_0 A)^{\langle s\rangle}=
\int_0C_1(s)=\tilde \gamma_1(s)
\in \mathcal{SD}$ (o\`u 
$\alpha^{\langle s\rangle}=\alpha\circ\alpha\circ\cdots\circ\alpha$ si
$s\in\mathbb N$
et o\`u
$\alpha^{\langle s\rangle}=
\alpha^{\langle-1\rangle}\circ\alpha^{\langle-1\rangle}\circ\cdots
\circ  \alpha^{\langle-1\rangle}$ si $s\in-\mathbb N$) 
pour tout $s\in\mathbb Z$.
En utilisant la structure de groupe de Lie sur $\mathcal{SI}$ (voir
les chapitres suivants), on peut
m\^eme consid\'erer $s$ \`a valeurs dans $\mathbb C$.
\end{rem}

\begin{rem} L'interpolation du corollaire \ref{corinterpol}
et la remarque pr\'ecedente
s'\'etendent sans difficult\'e aux groupes $
\mathcal U_{>0}={\mathbb R}_{>0}+\mathfrak
m
\subset \mathcal U$ et $\mathcal D_{>0}=x({\mathbb R}_{>0}+\mathfrak
m)\subset \mathcal D$. En travaillant dans le rev\^etement universel
$\tilde\mathcal I$ d\'ecrit dans la remarque \ref{remrvtmntuniv}, on
peut interpoler entre les groupes 
$\tilde \mathcal U$ et $\tilde \mathcal D$ obtenus \`a partir
de $\mathcal U$ et $\mathcal D$ en consid\'erant des rel\`evements
r\'eels des argument de  $A(0)$ et $\alpha'(0)$ pour $A\in
\mathcal U,\alpha\in\mathcal D$. 
\end{rem}

\begin{rem} 
Les bijections $(A,xA^{\tau_1})\longmapsto (A,xA^{\tau_2})$ et 
$(A,\int_0A^{\tau_1})\longmapsto
(A,\int_0A^{\tau_2})$ 
sous-jacentes \`a l'interpolation ne sont
pas des homomorphismes de groupes pour $\tau_1\not= \tau_2$.
\end{rem}

\begin{rem} Les deux fa\c cons d'interpoler entre 
le groupe commutatif $\mathcal{SU}=1+x{\mathbb C}[[x]]$ et 
le groupe non-commutatif $\mathcal{SD}=x+x^2{\mathbb C}[[x]]$
utilisent les deux bijections ensemblistes 
naturelles $A\longmapsto xA$  et $A\longmapsto \int_0A$
entre
$\mathcal{SU}$ et $\mathcal{SD}$. On peut \'evidemment 
les combiner avec des d\'eformations $\tau\longmapsto 
A_\tau$ de $A$ pour obtenir d'autres formules, par exemple
pour d'autres bijections entre $\mathcal{SU}$ et $\mathcal{SD}$. 
Pour l'interpolation entre $\frac{1}{A}$ et $\frac{(\int_0A)^{\langle
    -1\rangle}}{x}$, on peut ainsi par
exemple \'egalement prendre
$$\tau\longmapsto \frac{1}{A_\tau\circ(xA_\tau^\tau)^{\langle -1\rangle}}
$$
avec $A_\tau=\sum_{n=0}^\infty A_n\frac{x^n}{(n+1)^\tau}$.

De la m\^eme mani\`ere, la fonction
$$\tau\longmapsto \frac{1}{\tilde A_\tau\circ(\int_0\tilde
A_\tau^\tau)^{\langle -1\rangle}}
$$
avec $\tilde A_\tau=\sum_{n=0}^\infty (n+1)^\tau A_n x^n$, interpole entre
$\frac{1}{A}$ et $\frac{(xA)^{\langle
    -1\rangle}}{x}$.
\end{rem}

\begin{rem} Mentionnons qu'il existe beaucoup d'autres
formules d'interpolation qui sont cependant moins naturelles
car non reli\'ees \`a la structure de groupe de $\mathcal{SI}$.
Un exemple est donn\'e par la formule
$$\tau\longmapsto
\frac{1}{A\circ\left(x\left(\frac{1}{x}\int_0A\right)^\tau\right)}$$
qui interpole entre $\frac{1}{A}\in \mathcal{SU}$ et
$\frac{1}{x}\left(\int_0A\right)^{\langle -1\rangle}$.
\end{rem} 
\section{La repr\'esentation matricielle de ${\mathcal I}$}
\label{sectreprmatr}

Pour $(A,\alpha)\in \mathcal I=\mathcal{U}\rtimes \mathcal{D}$, notons
$\rho(A,\alpha)$ la matrice infinie de coefficients
$$\left(\rho(A,\alpha)\right)_{i,j}=[x^i]\left(A\alpha^j\right),\ 0\leq
i,j\in\mathbb N$$
o\`u $[x^i]\left(A\alpha^j\right)$ d\'esigne le coefficient
$\gamma_i$ de la s\'erie formelle 
$A\alpha^j=\sum_{n=j}^\infty \gamma_n x^n$.
Notons $\rho({\mathcal I})=\{\rho(A,\alpha)\ \vert\ (A,\alpha)\in 
\mathcal I\}$ 
l'ensemble de ces matrices. L'application
$(A,\alpha)\longmapsto \rho(A,\alpha)$ entre $\mathcal I$ et 
$\rho(\mathcal I)$ est bijective car $A=\sum_{n=0}^\infty
\left(\rho(A,\alpha)\right)_{n,0}x^n\in\mathcal U$ est la s\'erie 
g\'en\'eratrice de la colonne d'indice $0$ dans $\rho(A,\alpha)$ et 
$A\alpha=\left(\sum_{n=1}^\infty
\left(\rho(A,\alpha)\right)_{n,1}x^n\right)$ est
la s\'erie g\'en\'eratrice de la colonne d'indice $1$ dans
$\rho(A,\alpha)$. Remarquons aussi que $\rho({\mathcal I})$ est form\'e
de matrices qui sont triangulaires inf\'erieures et infinies.
Le sous-ensemble $\rho(\mathcal{SI})$ associ\'e au sous-groupe
$\mathcal{SI}\subset \mathcal I$ consiste en les matrices unipotentes
de $\rho(\mathcal I)$.

\begin{thm} \label{thmreplin}
L'ensemble $\rho({\mathcal I})$ est un groupe de matrices et
l'application $(A,\alpha)\longmapsto \rho(A,\alpha)$ est un
isomorphisme entre $\mathcal I={\mathcal U}\rtimes {\mathcal D}$ et
$\rho({\mathcal I})$. 
\end{thm}

\begin{rem} \label{remLaurent}
On pourrait \'egalement consid\'erer le produit
  semi-direct $$\mathbb C[[x]][x^{-1}]^*\rtimes \mathcal D$$
obtenu en rempla\c cant le groupe 
$\mathcal U$ par le groupe des \'el\'ements non-nuls
$\mathbb C[[x]][x^{-1}]^*$ du corps des s\'eries de Laurent.
Ce groupe admet une repr\'esentation matricielles dans les
 matrices biinfinies en associant \`a 
$(A,\alpha)\in \mathbb C[[x]][x^{-1}]^*\rtimes \mathcal D$
la matrice $\rho(A,\alpha)$ avec coefficients
$(\rho(A,\alpha))_{i,j}=[x^i](A\alpha^j)$ pour $i,j\in\mathbb Z$.
\end{rem} 

\begin{rem} L'homomorphisme
  $\varphi_{0,\lambda,0}:(A,\alpha)\longmapsto \left(A\left\lbrace
\frac{\alpha}{x}\right\rbrace^\lambda,\alpha\right)$ consid\'er\'e 
dans la proposition 
\ref{propvarphi} provient de la repr\'esentation lin\'eaire $\rho$ : 
En effet, pour $\lambda\in\mathbb N$ un entier naturel,
la matrice $\rho(\varphi_{0,\lambda,0}(A,\alpha))$ s'obtient en effa\c
cant les $\lambda$ premi\`eres lignes et colonnes de la matrice
triangulaire inf\'erieure $\rho(A,\alpha)$.
\end{rem}

{\bf Preuve du th\'eor\`eme \ref{thmreplin}} Le r\'esultat
d\'ecoule du calcul
$$\begin{array}{lcl}
(\rho(A,\alpha)\rho(B,\beta))_{i,j}&=&\sum_k[x^i](A\alpha^k)[x^k](B\beta^j)\\
&=&[x^i]\sum_kA(B\circ \alpha)(\beta\circ\alpha)^j\\
&=&
(\rho((A,\alpha)(B,\beta)))_{i,j}\end{array}$$
et du fait que $\rho(1,x)$ est la matrice identit\'e.
\hfill$\Box$

Dor\'enavant nous allons parfois utiliser la repr\'esentation 
fid\`ele $(A,\alpha)\longmapsto \rho(A,\alpha)$ pour identifier 
le groupe abstrait 
$\mathcal I$ avec sa repr\'esentation matricielle 
$\rho({\mathcal I})$.

\section{L'alg\`ebre de Lie de $\rho(\mathcal I)$}\label{sectalgLie}

Associons \`a une s\'erie formelle 
$\alpha=\sum_{n=0}^\infty \alpha_n x^n\in {\mathbb
  C}[[x]]$
les deux matrices triangulaires inf\'erieures infinies 
$$u_\alpha=\left(\begin{array}{ccccccccc}
\alpha_0\\
\alpha_1&\alpha_0\\
\alpha_2&\alpha_1&\alpha_0\\
\alpha_3&\alpha_2&\alpha_1&\alpha_0\\
\vdots&          &        &  &\ddots\end{array}\right),\quad
d_\alpha=\left(\begin{array}{ccccccccc}
0\\
0&\alpha_0\\
0&\alpha_1&2\alpha_0\\
0&\alpha_2&2\alpha_1&3\alpha_0\\
0&\alpha_3&2\alpha_2&3\alpha_1&4\alpha_0\\
\vdots&&&&&\ddots\end{array}\right)$$
dont les coefficients sont donn\'es par
$(u_\alpha)_{i,j}=\alpha_{i-j}$
et $(d_\alpha)_{i,j}=j\alpha_{i-j}$ pour $i,j\geq 0$. Remarquons qu'on
a $d_\alpha=u_\alpha d_1$ o\`u la matrice
$d_1$ associ\'ee \`a la s\'erie constante $1$
est la matrice diagonale infinie
de coefficients diagonaux $0,1,2,3,4,\dots$ les entiers naturels.
Notons 
$${\mathfrak i}=
\{u_\alpha+d_\beta\ \vert\ \alpha,\beta\in {\mathbb C}[[x]]\}$$
l'espace vectoriel engendr\'e par les matrices de la forme $u_\alpha,
d_\alpha$. Notons encore 
$$\mathfrak{si}\subset =
\{u_\alpha+d_\beta\ \vert\ \alpha,\beta\in \mathfrak m\}
\subset \mathfrak i$$
le sous-espace de codimension $2$ des matrices triangulaire
inf\'erieures strictes dans $\mathfrak i$.

\begin{thm} \label{algl} L'espace vectoriel ${\mathfrak i}$
est l'alg\`ebre de Lie de $\rho(\mathcal I)$. 
Le sous-espace $\mathfrak{si}\subset\mathfrak i$ 
correspond au sous-groupe $\rho(\mathcal{SI})$. Le crochet de Lie sur 
$\mathfrak i$ est donn\'e par les formules
$$\begin{array}{lcl}
\displaystyle [u_\alpha,u_\beta]=0\\
\displaystyle [d_\alpha,u_\beta]=u_{(x\alpha\beta')}\\
\displaystyle
[d_\alpha,d_\beta]=d_{(x(\alpha\beta'-\alpha'\beta))}\ .
\end{array}$$
et on a en particulier $[\mathfrak i,\mathfrak i]=\mathfrak{si}$.
\end{thm}

{\bf Preuve}  Pour $\alpha, \beta\in{\mathbb C}[[x]]$ 
et $h$ une variable  formelle consid\'erons
l'\'el\'ement de l'espace tangent $T_{(1,x)}(\mathcal I)$ en 
$\rho(1,x)\in\rho(\mathcal I)$ 
d\'efini par l'application $h\longmapsto (1+h\alpha,x(1+h\beta))$.
Un calcul \'el\'ementaire montre qu'on a 
$$\rho(1+h\alpha,x(1+h\beta))=\hbox{Id}+h(u_\alpha+d_\beta)+O(h^2)\ .$$
L'espace tangent en $\rho(1,x)\in\rho(\mathcal I)$
est donc donn\'e par l'espace vectoriel $\mathfrak i$
qui s'identifie alors \`a l'alg\`ebre de Lie de 
$\rho(\mathcal I)$. Comme le sous-groupe $\rho(\mathcal{SI})\subset
\rho(\mathcal I)$ est form\'e 
de matrices triangulaires inf\'erieures unipotentes, 
son alg\`ebre de Lie $\mathfrak{si}$ correspond aux \'el\'ements
triangulaires inf\'erieures strictes de $\mathfrak i$.

Le calcul du crochet de Lie r\'esulte des identit\'es
suivantes (qui ne font intervenir que des sommes finies):
$$
[u_\alpha,u_\beta]_{i,j}=\sum_k \alpha_{i-k}\beta_{k-j}-\beta_{i-k}
\alpha_{k-j}=0\ ,$$
$$\begin{array}{lcl}
[d_\alpha,u_\beta]_{i,j}&=&\sum_k k\alpha_{i-k}\beta_{k-j}-j
\sum_k\beta_{i-k}\alpha_{k-j}\\
&=&\sum_k k\alpha_{i-k}\beta_{k-j}-j
\sum_k\alpha_{i-k}\beta_{k-j}\\
&=&\sum_k (k-j)\alpha_{i-k}\beta_{k-j}\\
&=&\left(u_{(x\alpha\beta')}\right)_{i,j}\end{array}$$
et
$$\begin{array}{lcl}
[d_\alpha,d_\beta]_{i,j}&=&\sum_k
kj\left(\alpha_{i-k}\beta_{k-j}-\beta_{i-k}\alpha_{k-j}\right)\\
&=&j\sum_k (k-j)\left(\alpha_{i-k}\beta_{k-j}-\beta_{i-k}\alpha_{k-j}
\right)\\
&=&\left(d_{(x(\alpha\beta'-\alpha'\beta))}\right)_{i,j}\ .\end{array}$$
Ceci termine la preuve. \hfill$\Box$

Pour terminer ce chapitre, mentionnons sans donner les preuves 
faciles les faits suivants :

\begin{prop} \label{propfaitssuivants} (i)
Le sous-espace vectoriel $\mathfrak u=\{u_\alpha\ \vert\ \alpha
\in\mathbb C[[x]]\}\subset\mathfrak i$ est l'alg\`ebre de Lie 
(avec crochet nul) du sous-groupe commutatif 
$\rho(\mathcal U,x)\sim \mathcal U$ de $\rho(\mathcal I)$.

\ \ (ii) L'endomorphisme de groupe 
$\varphi_{\kappa,\lambda,\mu}(C,\gamma)=\left(C^\kappa\left\lbrace
\frac{\gamma}{x}\right\rbrace^\lambda\left\lbrace\gamma'\right
\rbrace^\mu,\gamma\right)$ de $\mathcal{SI}$ 
(voir proposition \ref{propvarphi})
correspond \`a (la restriction au sous-espace $\mathfrak{si}$ de)
l'endomorphisme d'alg\`ebre de Lie
(aussi not\'e) $\varphi_{\kappa,\lambda,\mu}
(u_\alpha+d_\beta)=u_{(\kappa\alpha+\lambda\beta+\mu(\beta+x\beta'))}
+d_\beta$ de
$\mathfrak i$ et provient d'un endomorphisme de groupe
$\tilde\varphi_{\kappa,\lambda,\mu}$ du ``rev\^etement universel'' $\tilde
\mathcal I$ de $\mathcal I$, obtenu en relevant les arguments
$\hbox{arg}(C(0)),\hbox{arg}(\gamma'(0))\in\mathbb R/(2\pi\mathbb Z)$
dans $\mathbb R$.

\ \ (iii) Le sous-espace vectoriel $\mathfrak d=\{d_\alpha\ \vert\ \alpha
\in\mathbb C[[x]]\}\subset\mathfrak i$ est l'alg\`ebre de Lie 
du sous-groupe 
$(1,\mathcal D)\sim \mathcal D$ de $\mathcal I$.
Plus g\'en\'eralement, les sous-espaces $\varphi_{0,\lambda,\mu}
(\mathfrak d)=
\{u_{(\lambda+\mu)\alpha+\mu x\alpha'}+d_\alpha\vert \alpha\in\mathbb
C[[x]]\}$ sont des sous-alg\`ebres
de Lie et les alg\`ebres de Lie
$\varphi_{0,\lambda,\mu}(\mathfrak d)\cap\mathfrak{si}$
correspondent aux sous-groupes de Lie
$\varphi_{0,\lambda,\mu}(1,\mathcal{SD})=\{\left(\left\lbrace\frac{\alpha}{x}\right\rbrace^\lambda\lbrace
  \alpha'\rbrace^\mu,\alpha\right)\ \vert\ \alpha\in \mathcal{SD}\}
\subset\mathcal{SI}$ qui sont tous isomorphes \`a $\mathcal{SD}$.

\ \ (iv) Les sous-groupes
$\mathcal{SG}(\tau)=\{(A,xA^\tau)\vert A\in\mathcal{SU}\}$ et
$\mathcal{SG}'(\tau)=\{(A,\int_0A^tau)\vert A\in\mathcal{SU}\}$
de $\mathcal{SI}$ utilis\'es dans le chapitre \ref{sectinterpol}
pour interpoler entre $\mathcal{SU}$ et $\mathcal{SD}$
correspondent aux sous-alg\`ebres de Lie
$\mathfrak{sg}(\tau)=\mathfrak g(\tau)\cap\mathfrak si$ 
et $\mathfrak{sg}'(\tau)=\mathfrak g(\tau)\cap\mathfrak si$ de
$\mathfrak{si}$ o\`u $\mathfrak g(\tau)=\{u_\alpha+\tau d_\alpha\vert
\alpha\in\mathbb C[[x]]\}$ et
$\mathfrak g'(\tau)=\{u_{(x\alpha)'}+\tau d_\alpha\vert
\alpha\in\mathbb C[[x]]\}$.
\end{prop}

\section{L'application exponentielle $\hbox{exp}:\mathfrak
  i\longrightarrow \rho(\mathcal I)$}\label{sectexp}

Comme $M=u_\alpha+d_\beta\in{\mathfrak i}$ est une matrice 
triangulaire inf\'erieure, l'exponentielle matricielle
$$\hbox{exp}(M)=\sum_{n=0}^\infty \frac{M^n}{n!}$$
converge vers une matrice triangulaire inf\'erieure inversible.
(Pour une matrice triangulaire inf\'erieure stricte
$M=u_\alpha+d_\beta\in\mathfrak{si}$, la situation est encore 
meilleure car le coefficient $(M^n)_{i,j}$ est nul pour $n>i-j$
et il n'y a plus besoin d'analyse.) 
Il r\'esulte de r\'esultats classiques que 
l'ensemble $\{\hbox{exp}(M)\ \vert\ M\in{\mathfrak i}\}$ 
s'identifie au groupe de Lie $\rho(\mathcal I)$ consid\'er\'e 
pr\'ec\'edemment.

Pour $\alpha,\beta\in{\mathbb C}[[x]]$, notons
$\hbox{Exp}(\alpha;\beta)=\sum_{n=0}^\infty M_{n,0}x^n
\in {\mathbb C}[[x]]$ la s\'erie
g\'en\'eratrice de la colonne d'indice $0$ dans la matrice
$M=\hbox{exp}(u_\alpha+d_\beta)\in \rho(\mathcal I)$.

\begin{prop} \label{propcoeffM}
Le coefficient $M_{i,j}$ (pour $0\leq i,j$) 
de la matrice $M=\hbox{exp}(u_\alpha+d_\beta)\in \rho(\mathcal I)$ est
donn\'e par la formule
$$M_{i,j}=[x^{i-j}]\hbox{Exp}(\alpha+j\beta;\beta)\ .$$
La s\'erie formelle 
$x^j\hbox{Exp}(\alpha+j\beta;\beta)=\sum_{i=0}^\infty 
M_{i,j}x^i$ est donc la s\'erie g\'en\'eratrice de la colonne d'indice
$j$ de $M$.
\end{prop}

\begin{cor} \label{corCgamma}
Si $\hbox{exp}(u_\alpha+d_\beta)=\rho(C,\gamma)\in \rho(\mathcal
  I)$ (pour $\alpha,\beta\in{\mathbb C}[[x],C\in\mathcal U,
\gamma\in\mathcal D$) alors 
$$C=\hbox{Exp}(\alpha;\beta)$$ et
$$C\gamma=x\ \hbox{Exp}(\alpha+\beta;\beta)\ .$$
\end{cor}

\begin{rem} On verra (voir l'assertion (iii) de la proposition 
\ref{idfondexp}) qu'on a \'egalement $\gamma=x\
\hbox{Exp}(\beta;\beta)$ si
$\hbox{exp}(u_\alpha+d_\beta)=\rho(C,\gamma)$.
\end{rem}

\begin{rem} \label{remexpcstes}
Comme $\hbox{exp}(u_\alpha)=\rho(e^\alpha,x)$, la fonction 
$\alpha\longmapsto \hbox{Exp}(\alpha;0)=e^\alpha$ est l'exponentielle
usuelle d'une s\'erie formelle $\alpha\in{\mathbb C}[[x]]$.

De mani\`ere similaire, si $\alpha$ et $\beta$ sont tous les deux des
fonctions constantes, la matrice $u_\alpha+d_\beta$ est diagonale
et on a donc $\hbox{Exp}(\alpha;\beta)=e^\alpha$ pour
$\alpha,\beta\in\mathbb C$.
\end{rem}

{\bf Preuve de la proposition \ref{propcoeffM}} La formule est vraie
pour $j=0$ par d\'efinition de la s\'erie formelle
$\hbox{Exp}(\alpha;\beta)$.

Pour $j>0$, on utilise les identit\'es
$$\begin{array}{lcl}
\left(\hbox{exp}(u_\alpha+d_\beta)\right)_{i,j}&=&
\left(\varphi_{0,j,0}(\hbox{exp}(u_\alpha+d_\beta))\right)_{i-j,0}\\
&=&\left(\hbox{exp}(\varphi_{0,j,0}(u_\alpha+d_\beta))\right)_{i-j,0}\\
&=&\left(\hbox{exp}(u_{(\alpha+j\beta)}+d_\beta)\right)_{i-j,0}\end{array}$$
o\`u $\varphi_{0,j,0}$ est l'homomorphisme de groupe (ou d'alg\`ebre)
de la proposition \ref{propvarphi}
qui consiste \`a effacer les $j$ premi\`eres lignes et colonnes
de la matrice triangulaire infinie
$\hbox{exp}(u_\alpha+d_\beta)\in\rho(\mathcal I)$.\hfill$\Box$

La preuve du corollaire \ref{corCgamma} est \'evidente. 

\begin{prop} \label{idfondexp}
(i) On a pour tout $\alpha,\beta\in\mathbb C[[x]]$
l'\'egalit\'e
$$\hbox{exp}(u_\alpha+d_\beta)=\rho(\hbox{Exp}(\alpha;\beta),
x\ \hbox{Exp}(\beta;\beta))\in \rho(\mathcal I)\ .$$

\ \  (ii) On a 
$$\alpha\longmapsto \hbox{Exp}(\alpha;0)=e^\alpha\ .$$

\ \ (iii) On a
$$\begin{array}{l}
\displaystyle 
\hbox{Exp}(\kappa\alpha+\lambda\beta+\mu(x\beta)';\beta)\\
\displaystyle \quad =
\hbox{Exp}(\alpha;\beta)^\kappa\hbox{Exp}(\beta;\beta)^\lambda
\left(\frac{d}{dx}(x\ \hbox{Exp}
(\beta;\beta))\right)^\mu\ .\end{array}$$

\ \ (iv) On a 
$$\hbox{Exp}((s+t)\alpha;(s+t)\beta)=\hbox{Exp}(s\alpha;s\beta)
\left(\hbox{Exp}
(t\alpha;t\beta)\circ(x\hbox{Exp}(s\beta;s\beta)\right)\ .$$

\ \ (v) On a 
$$\hbox{Exp}\left(\sum_{n=0}^\infty \alpha_n x^n;\beta\right)=
e^{\alpha_0}\ \hbox{Exp}\left
(\sum_{n=1}^\infty \alpha_n x^n;\beta\right)\ .$$

\ \ (vi) On a pour tout $K\in\mathbb C^*$
$$\hbox{Exp}(\alpha\circ(Kx);\beta\circ(Kx))
=\hbox{Exp}(\alpha;\beta)\circ(Kx)\ .$$
\end{prop}

\begin{rem} \`A cause de l'assertion (ii), on peut consid\'erer 
l'application 
$(\alpha,\beta)\longmapsto \hbox{Exp}(\alpha;\beta)$ 
comme une d\'eformation param\'etr\'ee par $\beta\in\mathbb C[[x]]$
de la fonction exponentielle usuelle
$$\alpha\longmapsto e^\alpha=\hbox{Exp}(\alpha;0)=\sum_{n=0}^\infty
\frac{\alpha^n}{n!}\ .$$

L'assertion (iv) se sp\'ecialise par exemple en 
$$x\hbox{Exp}(2\alpha;2\alpha)=(x\hbox{Exp}(\alpha;\alpha))\circ
(x\hbox{Exp}(\alpha;\alpha))$$
ce qui correspond bien \`a la loi de groupe sur le sous-groupe
$\varphi_{0,1,0}(\mathcal D)
=\{(A,xA)\ \vert\ A\in\mathcal U\}\sim \mathcal D$ 
de $\mathcal I$.

L'assertion (v) est triviale : Elle exprime seulement le fait 
que la matrice $\hbox{exp}(u_{\alpha_0})=e^{\alpha_0} \hbox{id}$ 
est centrale dans $\rho(\mathcal I)$.  

L'assertion (vi) (\'equivariance de $\hbox{Exp}$ par rapport aux
homoth\'eties inversibles $x\longmapsto Kx$ 
de $\mathbb C$) provient de l'automorphisme int\'erieur 
(de $\mathcal I$ et de $\mathfrak i$) donn\'e
en conjugant par la matrice diagonale
$\hbox{exp}(d_{\hbox{log}(K)})$, voir aussi remarque \ref{remautomtriv}.
\end{rem}

{\bf Preuve de la proposition \ref{idfondexp}} 
L'assertion (i) suit de la structure de produit semi-direct sur 
$\mathcal I$. On peut \'egalement s'en convaincre 
en consid\'erant l'homomorphisme $\varphi_{0,\lambda,0}
(C,\gamma)=(C\left\lbrace\frac{\gamma}{x}\right\rbrace^\lambda,\gamma)$
pour $\lambda\in \mathbb N$ grand. On obtient
ainsi l'identit\'e
$$C\left\lbrace\frac{\gamma}{x}\right\rbrace^\lambda=
\hbox{Exp}(\alpha+\lambda\beta;\beta)$$ si 
$\hbox{exp}(u_\alpha+d_\beta)=(C,\gamma)$. On a donc 
$$\gamma=x\ \hbox{lim}_{n\rightarrow\infty}\left(\hbox{Exp}
(\alpha+n\beta;\beta)\right)^{1/n}\ .$$

L'assertion (ii) est
l'observation que $\mathfrak u=\{u_\alpha\ \vert\ \alpha\in
\mathbb C[[x]]\}$ est l'alg\`ebre de Lie du groupe commutatif
$\mathcal U$.

L'assertion (iii) traduit le fait que $\varphi_{\kappa,\lambda,\mu}$
est un endomorphisme de l'alg\`ebre de Lie $\mathfrak i$,
voir la remarque \ref{remautomtriv}.

L'assertion (iv) suit de la loi de groupe sur $\mathcal I$.

L'assertion (v) est une cons\'equence du fait que $u_1$
est la matrice identit\'e.

L'assertion (vi) provient de la conjugaison par la matrice diagonale
de coefficients diagonaux $1,K,K^2,K^3,\dots$ une progression 
g\'eom\'etrique.\hfill$\Box$  

{\bf Exemple} Pour $\alpha=\beta=x\in {\mathbb C}[[x]]$, nous avons
$\hbox{Exp}(x;x)=\frac{1}{1-x}=\sum_{n=0}^\infty x^n$ car
$$\hbox{exp}(u_x+d_x)=\hbox{exp}\left(\begin{array}{ccccc}
0\\1&0\\
0&2&0\\
0&0&3&0\\
\vdots&&&\ddots\end{array}\right)=
\left(\begin{array}{ccccccc}
1\\
1&1\\
1&2&1\\
1&3&3&1\\
1&4&6&4&1\\
\vdots&&&&&\ddots\end{array}\right)$$
est la matrice triangulaire inf\'erieure unipotente associ\'ee au
triangle de Pascal.

\section{\'Equations diff\'erentielles}\label{sectequadiff}

\begin{prop} \label{diffeqpourExp} (i) Pour 
$\alpha=\sum_{n=N_\alpha}^\infty \alpha_nx^n,
\beta=\sum_{n=N_\beta}^\infty\beta_n x^n\in \mathbb C[[x]]$ avec
$\alpha_{N_\alpha},\beta_{N_\beta}\in \mathbb C^*$ et $N_\beta<\infty$,
les s\'eries formelles $y=\hbox{Exp}(\alpha;\beta),
z=\hbox{Exp}(\beta;\beta)\in\mathbb C[[x]]$ sont les 
solutions formelles uniques du syst\`eme d'\'equations diff\'erentielles
$$\left\lbrace\begin{array}{ll}
\displaystyle x\beta y'=y\{\alpha\circ (xz)-\alpha\}\\
\displaystyle
x\beta z'=z\{\beta\circ (xz)-\beta\}\end{array}\right.$$
qui v\'erifient
$$y=\left\lbrace\begin{array}{ll}
e^{\alpha_0}+\mathfrak m&\hbox{si }\alpha_0\not=0\\
1+\frac{e^{N_\alpha\beta_0}-1}{N_\alpha\beta_0}\alpha_{N_\alpha}x^{N_\alpha}+
\mathfrak m^{N_\alpha+1}&\hbox{sinon}\end{array}\right.$$
(avec la convention
$\frac{e^{N_\alpha\beta_0}-1}{N_\alpha\beta_0}=1$ si $N_\alpha\beta_0=0$)
et 
$$z=\left\lbrace\begin{array}{ll}
e^{\beta_0}+\mathfrak m&\hbox{si }\beta_0\not=0\\
1+\beta_{N_\beta}x^{N_\beta}+
\mathfrak m^{N_\beta+1}&\hbox{sinon}\end{array}\right.$$

\ \ (ii) Les fonctions 
$$\begin{array}{lcl}y_s&=&\hbox{Exp}(s\alpha;s\beta)\\
y_s'&=&\frac{d}{dx}\hbox{Exp}(s\alpha;s\beta)\\
z_s&=&\hbox{Exp}(s\beta;s\beta)\\
z_s'&=&\frac{d}{dx}\hbox{Exp}(s\beta;s\beta)\end{array}$$
sont solution du syst\`eme d'\'equations diff\'erentielles
$$\left\lbrace\begin{array}{lcl}
\frac{\partial}{\partial s}y_s&=&\alpha y_s+x\alpha y_s'=y_s\alpha
\circ(xz_s)\\
\frac{\partial}{\partial s}z_s&=&\beta z_s+x\beta z_s'=z_s\beta
\circ(xz_s)\ .\end{array}\right.$$
\end{prop}

\begin{rem} Le param\`etre $s$ se simplifie 
dans le syst\`eme d'\'equations
diff\'erentielles 
$$\left\lbrace\begin{array}{ll}
\displaystyle x\beta y'=y\{\alpha\circ (xz)-\alpha\}\\
\displaystyle
x\beta z'=z\{\beta\circ (xz)-\beta\}\end{array}\right.$$
pour $y=\hbox{Exp}(s\alpha;s\beta)$ et 
$z=\hbox{Exp}(s\beta;s\beta)$. Le param\`etre $s$ 
n'intervient donc que dans les conditions initiales.

Le syst\`eme d'\'equations diff\'erentielles
pour $y_s=\hbox{Exp}(s\alpha;s\beta)$ et $z_s=
\hbox{Exp}(s\beta;s\beta)$ est autonome en $s$ car il provient d'un
flot sur le groupe de Lie $\rho(\mathcal I)$.

Sous des hypoth\`eses de convergence convenables, le calcul
d'une \'evaluation en $x=x_0\in \mathbb C$ de
$y=\hbox{Exp}(\alpha;\beta)$
peut \'egalement se faire en d\'eterminant la valeur en $s=1$ 
de la fonction $u=u(s)$ pour $u,w$ les fonctions (correspondantes
\`a l'\'evaluation en $x=x_0$ de $\hbox{Exp}(s\alpha;s\beta)$
et  $\hbox{Exp}(s\alpha;s\beta)$)
d\'etermin\'ees par le probl\`eme de Cauchy
$$\left\lbrace \begin{array}{l}
u'=u\ \alpha\circ(x_0w)\\
w'=w\ \beta\circ(x_0w)\end{array}\right.$$
avec conditions initiales $u(0)=w(0)=1$. 
\end{rem}

{\bf Exemple} (i) $\hbox{Exp}(sx;s)=e^{(e^s-1)x}=
1+(e^s-1)x+\frac{(e^s-1)^2}{2}2x^2+\dots$. En effet,
posons $y=e^{(e^s-1)x}=1+(e^s-1)x+\dots$ et $z=e^s$. Les fonctions 
$y$ et $z$ v\'erifient alors le syst\`eme d'equations
diff\'erentielles
$$\left\lbrace\begin{array}{l}xsy'=y\{sxe^s-sx\}\\
xsz'=z\{s-s\}=0\end{array}\right.$$
(o\`u l'on d\'erive par rapport \`a la variable $x$)
correspondant \`a $\alpha=sx$ et $\beta=s$.
De plus, on a $N_\alpha=1,\alpha_1=s,N_\beta=0,\beta_0=s$
et $y,z$ sont donc bien de la forme
$y\in 1+\frac{e^{1\cdot s}-1}{1\cdot s}sx+\mathfrak m^2$,
$z\in e^s+\mathfrak m$, voir aussi la remarque \ref{remexpcstes}
pour l'\'egalit\'e $z=e^s=\hbox{Exp}(s;s)$.

\ \ (ii) La fonction 
$y=\hbox{Exp}\left(\frac{x}{1-x},\frac{x}{1-x}\right)$
est l'unique solution de la forme $y\in 1+x+\mathfrak m^2$
de l'\'equation diff\'erentielle $$
x\frac{x}{1-x}y'=y\left\lbrace\frac{xy}{1-xy}-\frac{x}{1-x}\right
\rbrace$$
(correspondante \`a $\beta=x/(1-x)$) qui se simplifie en
$$x\lbrace 1-xy\rbrace y'=y\lbrace y-1\rbrace\ .$$
La suite $c_0=y(0),c_1=y'(0),c_2=y''(0),\dots,c_n=y^{(n)}(0),\dots$ 
des d\'eriv\'ees en $0$ de la solution
$y(x)=\hbox{Exp}\left(\frac{x}{1-x},\frac{x}{1-x}\right)=
\sum_{n=0}^\infty c_n\frac{x^n}{n!}$
commence par
$$1,1,4,27,260,3270,50904,946134,20462896,\dots\ .$$
 
La suite $\tilde c_1=-\left(\frac{1}{y}\right)',
\tilde c_2=-\left(\frac{1}{y}\right)'',\dots,\tilde c_n=
-\left(\frac{1}{y}\right)^{(n)},\dots$
associ\'ee aux d\'eriv\'ees d'ordre $\geq 1$ de
$\frac{1}{y}=1-\sum_{n=1}^\infty \tilde c_n\frac{x^n}{n!}$ 
admet une interpr\'etation combinatoire et
correspond \`a la suite A38037 dans \cite{EIS}.

{\bf Preuve de la proposition \ref{diffeqpourExp}} 
L'assertion (iv) de la proposition \ref{idfondexp} donne
$$\frac{\partial}{\partial
  s}\hbox{Exp}((s+t)\alpha;(s+t)\beta)\vert_{s=0}
=\alpha\
\hbox{Exp}(t\alpha;t\beta)+\left(\frac{d}{dx}\hbox{Exp}
(t\alpha;t\beta)\right)\lbrace x\beta\rbrace\ .$$
On trouve de m\^eme
$$\frac{\partial}{\partial
  t}\hbox{Exp}((s+t)\alpha;(s+t)\beta)\vert_{t=0}
=\hbox{Exp}(s\alpha;s\beta)\alpha\circ\left(x\hbox{Exp}(s\beta;s\beta)
\right)\ .$$
En posant $t=s=1,y=\hbox{Exp}(\alpha;\beta),
z=\hbox{Exp}(\beta;\beta)$ et en \'egalant les
deux identit\'es, on trouve la premi\`ere 
\'equation du syst\`eme diff\'erentiel. La preuve pour la 
deuxi\`eme \'equation est analogue.

La forme des solutions vient de la d\'efinition de
$y$, respectivement de $z$, comme s\'erie g\'en\'eratrice de la
premi\`ere colonne de $\hbox{exp}(u_\alpha+d_\beta)$,
respectivement de $\hbox{exp}(u_\beta+d_\beta)$,
et du petit calcul
$$\hbox{exp}\left(\begin{array}{cc}0&0\\
\alpha_{N_\alpha}&N_\alpha \beta_0\end{array}\right)=
\left(\begin{array}{cc}1&0\\
\alpha_{N_\alpha}\frac{e^{N_\alpha \beta_0}-1}{N_\alpha\beta_0}  
&e^{N_\alpha \beta_0}\end{array}\right)\ .$$
Ceci termine la preuve de l'assertion (i). 

L'assertion (ii) est une cons\'equence
facile de ce qui pr\'ec\`ede.\hfill$\Box$

\begin{rem} Le syst\`eme diff\'erentiel de la proposition
 \ref{diffeqpourExp}
d\'eg\'en\`ere pour $\beta=0$. Ceci est d\^u
au fait qu'il a \'et\'e d\'eriv\'e par l'utilisation
de la non-commutativit\'e du groupe ${\mathcal I}$.
\end{rem} 

\subsection{La fonction $x\ \hbox{Exp}(\beta;\beta)$}\label{expbetabeta}

L'application $\beta\longmapsto x\ \hbox{Exp}(\beta;\beta)$
correspond \`a l'exponentielle entre la sous-alg\`ebre
de Lie $\mathfrak d\subset \mathfrak i$ et le groupe de Lie
$\mathcal D$, identifi\'e au sous-groupe
form\'e des \'el\'ements $(1,\alpha),\ \alpha\in\mathcal D$ 
dans $\mathcal I$.

\begin{prop} \label{propeqdiffbetabeta}
Soit $\beta\in \mathbb C[[x]]$ une s\'erie formelle
non-nulle pour laquelle la s\'erie formelle de Laurent 
$\frac{1}{x\beta}$ est sans r\'esidu pour le p\^ole \`a
l'origine. Alors
$$F(x\ \hbox{Exp}(s\beta;s\beta))=F(x)+s$$
pour $F$ une primitive de $\frac{1}{x\beta}$.

Ce r\'esultat reste valable pour une s\'erie formelle 
$\beta\in\mathbb C[[x]]$ quelconque pour une d\'etermination
convenable du logarithme apr\`es simplification
de la singularit\'e essentielle due au r\'esidu.
\end{prop}

\begin{cor} \label{expalg}
Si $\frac{1}{x\beta}$ admet une primitive m\'eromorphe
alg\'ebrique pour $\beta\in\mathfrak m$, alors
$\hbox{Exp}(s\beta;s\beta)$ est alg\'ebrique.

En particulier, $\hbox{Exp}(s\beta;s\beta)$ est alg\'ebrique
pour $\beta\in x\mathbb C[[x]]$ une fraction rationelle telle
que $\frac{1}{x\beta}$ est sans r\'esidu en tout 
point de $\mathbb C$. 
\end{cor}

\begin{rem} La s\'erie de Laurent $\frac{1}{x\beta}$ poss\`ede un 
p\^ole d'ordre $1$ \`a l'origine si $\beta=\beta_0+\beta_1x+\dots$ 
avec $\beta_0\not=0$. La partie singuli\`ere de $F$
provient alors de $\frac{1}{\beta_0}\hbox{log}(x)$ et on a
$$\begin{array}{l}
\displaystyle 
\frac{1}{\beta_0}\hbox{log}(x\ \hbox{Exp}(s\beta;s\beta))-
\frac{1}{\beta_0}\hbox{log}(x)=\frac{1}{\beta_0}\hbox{log}(\hbox{Exp}(s\beta;s\beta))\\
\displaystyle 
=\frac{1}{\beta_0}\hbox{log}(e^{s\beta_0}(1+\dots))=s+\hbox{log}
(1+\dots)\in s+\mathfrak m\end{array}$$
qui est une s\'erie formelle ordinaire en $0$.

Dans le cas d'une singularit\'e $r \hbox{log}(x)$ de $F$ provenant
du r\'esidu d'un p\^ole multiple \`a l'origine, on a
$\beta\in\mathfrak m$ et le calcul
$$\begin{array}{l}
\displaystyle 
r\ \hbox{log}(x\ \hbox{Exp}(s\beta;s\beta))-
r\ \hbox{log}(x)=r\ \hbox{log}(\hbox{Exp}(s\beta;s\beta))\\
\displaystyle 
=r\ \hbox{log}(1+b_1x+\dots)\in \mathfrak m\end{array}$$
montre \'egalement que les singularit\'es essentielles 
des deux c\^ot\'es se simplifient dans la proposition 
\ref{propeqdiffbetabeta}
\end{rem}

{\bf Preuve de la proposition \ref{propeqdiffbetabeta}} 
Le changement de 
variable $Z=xz=x\ \hbox{Exp}(\beta;\beta)$ transforme 
l'\'equation diff\'erentielle $x\beta
z'=z\lbrace\beta\circ(xz)-\beta\rbrace$ de la proposition 
\ref{diffeqpourExp} en
$$xZ'\beta(x)=Z\beta(Z)$$
qui est \`a variables s\'epar\'ees. Apr\`es int\'egration, on
obtient
$$F(x\ \hbox{Exp}(s\beta;s\beta))=F(x)+h(s)$$
pour $F$ tel que $F'=\frac{1}{x\beta}$ et
o\`u $h(s)$ ne d\'epend que de $s$.
Dans le cas o\`u $\frac{1}{x\beta}$ est sans r\'esidu \`a l'origine,
on a 
$$F(x)=-\frac{1}{N\beta_Nx^N}+\dots\in\frac{1}{x^N}\mathbb C[[x]]$$
si $\beta=\sum_{n=N}^\infty\beta_nx^n$ avec $\beta_N\not=0$ 
pour un certain entier $N\geq 1$.
En utilisant la condition initiale $x\ \hbox{Exp}(s\beta;s\beta)=
x+s\beta_Nx^{N+1}+\dots$ on obtient
$$[x^0]F(x\ \hbox{Exp}(s\beta;s\beta))=[x^0]\left(
-\frac{1}{N\beta_Nx^N}\left(1-s\beta_Nx^N\right)^N\right)=s$$
ce qui montre $h(s)=s$ et
d\'emontre le r\'esultat dans le cas o\`u $\frac{1}{x\beta}$
est sans r\'esidu.

Dans le cas $\beta\in\mathfrak m$ et
$\frac{1}{x\beta}$ a un r\'esidu non-nul, la fonction $h(s)$
peut s'\'evaluer en calculant la limite $x\rightarrow 0$
et le r\'esultat suit 
de l'\'egalit\'e
$$\hbox{log}(x\ \hbox{Exp}(s\beta;s\beta))=\hbox{log}(x)+
\hbox{log}(\hbox{Exp}(s\beta;s\beta))=\hbox{log}(x)+
\hbox{log}(1+s\beta_1x+\dots)=\ .$$

Dans le cas o\`u $\beta=\beta_0+\dots$ avec $\beta_0\not=0$,
on a $F=\frac{1}{\beta_0}\hbox{log}(x)+\mathfrak m$ et on obtient 
encore $h(s)=s$.\hfill$\Box$

{\bf Preuve du corollaire \ref{expalg}}
La preuve de la premi\`ere partie est imm\'ediate.
La deuxi\`eme partie suit de l'observation qu'une primitive d'une
telle fraction rationelle est une fraction rationelle.
\hfill$\Box$

{\bf Exemples} 
(i) Cet exemple est inspir\'e du chapitre 3 de \cite{CN} 
o\`u il est consid\'er\'e \`a cause de son int\'er\^et pour la
biologie mol\'eculaire, voir les r\'ef\'erences
indiqu\'ees dans \cite{EIS} pour la suite A4148. On a
$$\begin{array}{l}
\displaystyle
\hbox{Exp}\left(s\frac{x}{1-x^2};s\frac{s}{1-x^2}\right)\\
\displaystyle \qquad =
\frac{1-sx+x^2-\sqrt{1-2sx+(s^2-2)x^2-2sx^3+x^4}}{2x^2}\ .\end{array} $$
En effet, il suffit
de v\'erifier l'\'egalit\'e de la proposition \ref{propeqdiffbetabeta} 
avec $F=-\frac{1}{x}-x$ une primitive de 
$\frac{1}{x}\left(\frac{x}{1-x^2}\right)^{-1}=\frac{1-x^2}{x^2}$.

Pour $s=1$, on obtient la s\'erie
$$1+x+x^2+2x^3+4x^4+8x^5+17x^6+37x^7+82x^8+$$
dont la suite des coefficients est A4148 de \cite{EIS}.

\ \ (ii) Plus g\'en\'eralement, les fractions rationnelles 
$\beta=\frac{x(1+x^2)^k}{1-x^{2(k+1)}}$ semblent satisfaire les
conditions du corollaire \ref{expalg} pour tout $k\in\mathbb N$.

Les cas $k=0$ et $k=1$ redonnent l'exemple ci-dessus. 

Pour $k=2,3$ la fonction
$$Z=x\ \hbox{Exp}\left(s\frac{x(1+x^2)^k}{1-x^{2(k+1)}};
s\frac{x(1+x^2)^k}{1-x^{2(k+1)}}\right)$$
satisfait l'\'equation alg\'ebrique
$$\begin{array}{l}
\displaystyle 
x(1+x^2)(1+(k+1)Z^2+Z^4)\\
\displaystyle \qquad -(1-sx+(k+1)x^2-sx^3+x^4)Z(1+Z^2)\ .\end{array}$$

Le cas $k=2,\ s=1$ donne
$$Z=x+x^2+x^3+3x^4+7x^5+14x^6+33x^7+81x^8+\dots$$
et correspond \`a la fonction alg\'ebrique 
$$\left(\frac{1-\sqrt{1-4x^2}}{2x}\right)\circ
\left(\frac{(1-x+x^2)-\sqrt{1-2x-x^2-2x^3+x^4}}{2x}\right)\circ
\left(\frac{x}{1+x^2}\right)\ .$$

Le d\'eveloppement en s\'erie de $Z$ pour $k=3,\ s=1$ commence par 
$$Z=x+x^2+x^3+4x^4+10x^5+22x^6+61x^7+165x^8+\dots\ .$$

\section{Un d\'eveloppement en s\'erie pour
$\hbox{Exp}(s\alpha;s\beta)$}\label{sectdevserie}

Pour $\alpha=\sum_{n=0}^\infty \alpha_n x^n,\beta=\sum_{n=0}^\infty 
\beta_n x^n\in {\mathbb
  C}[[x]]$ deux s\'eries formelles, posons $g_0=1$ 
et d\'efinissons ensuite les s\'eries
$g_1,g_2,\dots\in {\mathbb C}[[x]]$ (qui d\'ependent de 
$\alpha$ et de $\beta$)
r\'ecursivement par la formule
$$g_{n+1}=\alpha g_n+x \beta g'_n
\in {\mathbb
  C}[[x]]$$
o\`u $g'_n\in\mathbb C[[x]]$ est la s\'erie d\'eriv\'ee de $g_n$.

\begin{prop} \label{propformExp}
On a pour tout $s\in {\mathbb C}^*$ l'\'egalit\'e
$$\hbox{Exp}(s\alpha;s\beta)=\sum_{n=0}^\infty
g_n\frac{s^n}{n!}\ .$$
\end{prop}

\begin{rem} \label{formopdiff}
La proposition \ref{propformExp} peut \'egalement
s'\'enoncer sous la forme
$$\hbox{Exp}(s\alpha;s\beta)=\left(\sum_{n=0}^\infty
\frac{s^n}{n!}\left(\alpha+x\beta\frac{\partial}{\partial
    x}\right)\right)1\ .$$
\end{rem}

\begin{rem} Le coefficient $[x^N]\hbox{Exp}(\alpha;\beta)$
ne d\'epend que de $\alpha_0,\alpha_1,\dots,\alpha_N$ et de
$\beta_0,\beta_1,\dots,\beta_{N-1}$. De plus, pour $\alpha,\beta
\in \mathfrak m=x\mathbb C[[x]]$, cette d\'ependance est polynomiale
car $g_n\in\mathfrak m^n$.
\end{rem} 

\begin{rem} \label{pasdedivpourcoeffdeExp}
La formule $\hbox{Exp}(s\alpha;s\beta)=\sum_{n=0}^\infty
g_n\frac{s^n}{n!}$ ne pose jamais de probl\`eme de divergences pour le
calcul des coefficients $[x^N]\hbox{Exp}(\alpha;\beta)$.
En effet, la majoration facile
$$\vert[x^N]g_{n+1}\vert\leq\left(\sum_{k=0}^N\vert
  [x^k]g_n\vert\right)
\left(\sum_{k=0}^N\vert \alpha_k\vert+N\sum_{k=0}^{N-1}
\vert \beta_k\vert\right)$$
implique les majorations 
$$\vert [x^N]\hbox{Exp}(s\alpha;s\beta)\vert\leq
\sum_{n=0}^\infty \vert [x^N]g_n\vert\frac{\vert s\vert^n}{n!}
\leq e^{A_N\vert s\vert}$$
pour $A_N=\sum_{k=0}^N\vert \alpha_k\vert+N\sum_{k=0}^{N-1}
\vert \beta_k\vert$.
\end{rem}

\begin{rem} Les formules pour le calcul de
  $\hbox{Exp}(s\alpha;s\beta)$ donn\'ees par la proposition 
\ref{propformExp} sont particuli\`erement jolies dans le cas
$$\hbox{Exp}(s(x\beta)';s\beta)=
\hbox{Exp}(s(\beta+x\beta');s\beta)=\frac{d}{dx}\left(x\
  \hbox{Exp}(s\beta;s\beta)\right)$$
(voir l'assertion (iii) de la proposition \ref{idfondexp} pour 
la derni\`ere \'egalit\'e) o\`u l'on obtient
$\hbox{Exp}(s(x\beta)';s\beta)=\left(\sum_{n=0}^\infty\frac{s^n}{n!}
\left(\frac{\partial}{\partial x}x\beta\right)\right)1=
\sum_{n=0}^\infty
g_n\frac{s^n}{n!}$ avec $g_0=1$ et
$g_{n+1}=\beta g_n+x\beta'g_n+x\beta g_n'=(x\beta g_n)'$.
\end{rem}

{\bf Preuve de la proposition \ref{propformExp}} 
La d\'emonstration consiste \`a calculer la s\'erie
g\'en\'eratrice de la premi\`ere
colonne dans la matrice 
$$\hbox{exp}(su_\alpha+sd_\beta)=\sum_{n=0}^\infty
(u_\alpha+d_\beta)^n\frac{s^n}{n!}\ .$$
Montrons par r\'ecurrence sur $n$ que le coefficient de 
$s^n$ dans cette s\'erie 
est donn\'ee par $\frac{g_n}{n!}$. 
Ceci est trivialement vrai pour $n=0$
car $(u_\alpha+d_\beta)^0$ est la matrice identit\'e par convention.
D\'esignons par ${\mathcal C}$ l'espace vectoriel des matrices
triangulaires inf\'erieures dont la premi\`ere colonne
est identiquement z\'ero. Il suffit alors de d\'emontrer
l'identit\'e
$$(u_\alpha+d_\beta)u_{g_n}=u_{g_{n+1}}\pmod{\mathcal C}\ .$$
Or elle r\'esulte du calcul
$$\begin{array}{lcl}
(u_\alpha+d_\beta)u_{g_n}&\equiv
&u_\alpha u_{g_n}+d_\beta u_{g_n}-u_{g_n}d_\beta\pmod{\mathcal C}\\
&=&u_{\alpha g_n}+[d_\beta,u_{g_n}]\\
&=&u_{\alpha g_n}+u_{x\beta g_n'}=u_{g_{n+1}}\end{array}$$
car $u_{g_n}d_\beta\in\mathcal C$ et $[d_\beta,u_{g_n}]= 
u_{x\beta g_n'}$ par les formules du th\'eor\`eme \ref{algl}.\hfill$\Box$

{\bf Exemple} Pour $\alpha=\beta=x$, on obtient facilement
$g_n=n!x^n$ par r\'ecurrence sur $n$.
On a donc $\hbox{Exp}(sx;sx)=\sum_{n=0}^\infty n!x^n\frac{s^n}{n!}=
\frac{1}{1-sx}$ ou encore $\hbox{Exp}(s\kappa
x;sx)=\left(\frac{1}{1-sx}\right)^\kappa$ en utilisant l'assertion
(iii) de la proposition \ref{idfondexp}. Plus g\'en\'eralement,
on a $\hbox{Exp}(s\kappa x^a;sx^a)=\left(\frac{1}{1-sax^a}
\right)^{\kappa/a}$ pour $a\in\{1,2,3,\dots\}$.

\section{Fonctions r\'eciproques}\label{sectfctrecipr}

Consid\'erons les polyn\^omes
$P_1=1,P_2=1+x,P_3=1+\frac{3x}{2}+\frac{3x^2}{2},\dots,
P_n=\sum_{k=1}^n \frac{(nx)^{k-1}}{k!},\dots\subset\mathbb Q[x]$.
Notons $\mathcal R_n=\{\xi\in\mathbb C\ \vert\ P_n(\xi)=0\}$
l'ensemble des racines de $P_n$ et consid\'erons la r\'eunion
d\'enombrable $\mathcal R=\cup_{n=1}^\infty \mathcal R_n\subset
\mathbb C$. Remarquons l'in\'egalit\'e
$\frac{(2n)^k}{(k+1)!}\geq 2 \frac{(2n)^{k-1}}{k!}$ 
pour $k<n$, qui implique l'inclusion 
$\mathcal R\subset\{z\in\mathbb C\ \vert\ \parallel z\parallel<2\}$.

\begin{thm} \label{thmrecipr} (i) Pour $\beta\in\mathbb C[[x]]$ 
avec $\beta_0\not\in \mathcal R$ et 
$\gamma=\sum_{n=0}^\infty \gamma_nx^n\in\mathcal U=\mathbb
C^*+\mathfrak m$ donn\'es, il existe 
$\alpha\in\mathbb C[[x]]$ tel que $\hbox{Exp}(\alpha;\beta)=\gamma$.

De plus, si $\hbox{Exp}(\alpha;\beta)=\hbox{Exp}(\tilde
\alpha;\beta)$, alors il existe $k\in\mathbb Z$ tel que 
$\tilde \alpha=\alpha+2ik\pi$.

\ \ (ii) Pour $
\gamma\in\mathcal U=\mathbb C^*+\mathfrak m$ fix\'e, il existe 
$\beta\in\mathbb C[[x]]$ tel que $\hbox{Exp}(\beta;\beta)=\gamma$.

De plus, si $\hbox{Exp}(\beta;\beta)=\hbox{Exp}(\tilde
\beta;\tilde \beta)$ avec $\beta=\sum_{j=0}^\infty \beta_nx^n,
\tilde \beta=\sum_{n=0}^\infty \tilde \beta_nx^n$, 
alors il existe $k\in\mathbb Z$ tel que 
$\tilde \beta_0=\beta_0+2ik\pi$.
\end{thm}

\begin{cor} \label{corexpsurj}
L'application de $\mathbb C[[x]]\times \mathbb C[[x]]$
dans $\mathcal I$
d\'efinie par 
$$(\alpha,\beta)\longmapsto
(\hbox{Exp}(\alpha;\beta),x\hbox{Exp}(\beta;\beta))$$
est surjective.

Sa restriction \`a $\mathfrak m\times \mathfrak m$ est une bijection
entre $\mathfrak m\times \mathfrak m$ et $\mathcal{SI}$.
\end{cor}

\begin{rem} (i) (cf. remarque \ref{remexpcstes})
Pour $\beta=0$ et $\gamma=\sum_{n=0}^\infty \gamma_nx^n\in\mathcal U$, 
une s\'erie 
$\alpha\in\mathbb C[[x]]$ telle que $\hbox{Exp}(\alpha;0)=\gamma$ est
de la forme $\hbox{log}(\gamma_0)\hbox{log}\left(
\frac{1}{\gamma_0}\gamma\right)$
o\`u $\hbox{log}(\gamma_0)\in\mathbb C$ est une d\'etermination 
du logarithme et o\`u $\hbox{log}\left(\frac{1}{\gamma_0}\gamma
\right)\in \mathfrak m$
est la d\'etermination principale du logarithme de
$\frac{1}{\gamma_0}\gamma\in \mathcal{SU}=1+\mathfrak m$.

\ \ (ii) Si $\hbox{Exp}(\beta;\beta)=\hbox{Exp}(\tilde
\beta;\tilde \beta)$ avec $\beta=\sum_{j=0}^\infty \beta_nx^n,
\tilde \beta=\sum_{n=0}^\infty \tilde \beta_nx^n$, alors 
les diff\'erences $\tilde \beta_n-\tilde \beta_n$ ne
sont g\'en\'eralement pas nulles pour $n\geq 1$.
\end{rem} 

{\bf Preuve du th\'eor\`eme \ref{thmrecipr}} 
Prouvons l'assertion (i) en supposant d'abord 
$\gamma\in 1+\mathfrak m$. La suite $g_0,g_1,\dots\in\mathbb C[[x]]$
d\'efinie par $g_0=1$ et $g_{n+1}=\alpha g_n
+x\beta g_n'$ \`a partir de $\alpha=\sum_{n=1}^\infty \alpha_nx^n
\in\mathfrak m$ v\'erifie 
$g_n\in\mathfrak m^n$. Pour $\beta$
fix\'e, le coefficient $[x^n]\hbox{Exp}(\alpha;\beta)$
est alors un polyn\^ome en $\alpha_1,\dots,\alpha_n$. 
Le coefficient $\alpha_n$ de $\alpha$ n'intervient que
de fa\c con lin\'eaire dans $[x^n]\hbox{Exp}(\alpha;\beta)$
et un petit calcul montre que sa contribution est donn\'ee par
$$\alpha_n\left(\frac{1}{1!}+
\frac{n\beta_0}{2!}+\dots+\frac{(n\beta_0)^{n-1}}{n!}
\right)=\alpha_nP_n(\beta_0)\ .$$
Comme $P_n(\beta_0)\not=0$ pour tout $n$, cette \'equation lin\'eaire
se r\'esoud toujours et on peut donc d\'eterminer 
r\'ecursivement les coefficients
$\alpha_1,\alpha_2,\dots$ de mani\`ere \`a avoir 
$\hbox{Exp}(\alpha;\beta)=\gamma$.

L'assertion (v) de la proposition \ref{idfondexp}
permet de ramener le cas g\'en\'eral $\gamma=\sum_{n=0}^\infty
\gamma_n x^n \in\mathcal U$ 
au cas particulier $\frac{1}{\gamma_0}\gamma\in
\mathcal{SU}=1+\mathfrak m$
d\'ej\`a trait\'e. Ceci termine la preuve de l'assertion (i).

Pour d\'emontrer l'assertion (ii), on commence par choisir $\beta_0
\in\mathbb C\setminus\mathcal R$ tel que $e^{\beta_0}=\gamma_0$.
Ceci est toujours possible car l'ensemble $\mathcal R$ est born\'e.
On proc\`ede ensuite comme ci-dessus. Plus pr\'ecis\'ement,
on r\'esoud l'\'equation
$$\hbox{Exp}(\beta-\beta_0;\beta)=e^{-\beta_0}\gamma$$en remarquant 
que la contribution de $\beta_n$ au polyn\^ome
$[x^n]\hbox{Exp}(\beta-\beta_0;\beta)$ en $\beta_0,\dots,\beta_n$
est encore donn\'e par $\beta_nP_n(\beta_0)$ pour $n\geq 1$.
\hfill $\Box$

{\bf Preuve du corollaire \ref{corexpsurj}} Soit $(C,\gamma)\in\mathcal I$
avec $C\in\mathcal U$ et $\gamma\in\mathcal D$. 
L'assertion (ii) du th\'eor\`eme \ref{thmrecipr} permet 
de trouver $\beta\in\mathbb C[[x]]$ avec $\beta_0\not\in\mathcal R$
et $x\ \hbox{Exp}(\beta;\beta)=\gamma$. On utilise ensuite
l'assertion (i) pour d\'eterminer $\alpha\in\mathbb C[[x]]$ tel que
$\hbox{Exp}(\alpha;\beta)=C$. Pour $(C,\gamma)\in\mathcal{SI}$
on peut prendre l'unique solution  $(\alpha,\beta)$ dans $\mathfrak m
\times \mathfrak m$.
\hfill$\Box$

\begin{rem} La preuve du th\'eor\`eme \ref{thmrecipr} fournit un
  algorithme de calcul pour $\alpha$ (respectivement $\beta$)
tel que $\hbox{Exp}(\alpha;\beta)=\gamma$ (respectivement 
$\hbox{Exp}(\beta;\beta)=\gamma$) pour $\beta,\gamma$ (respectivement
$\gamma$) convenable. En effet, soit $\tilde \alpha\in\mathbb C[[x]]$
une s\'erie telle que $e^{\tilde \alpha_0}=\gamma_0$ et 
$\gamma-\hbox{Exp}(\tilde\alpha;\beta)\equiv E_nx^n\pmod {\mathfrak
  m^{n+1}}$ pour un entier $n\geq 1$. Alors la s\'erie 
$$\overline \alpha=\tilde
\alpha+\frac{1}{\gamma_0}
E_n\left(\sum_{k=1}^{n}\frac{(n\beta_0)^{k-1}}{k!}\right)^{-1}
x^n$$
v\'erifie $\gamma-\hbox{Exp}(\overline \alpha;\beta)\equiv 0
\pmod{{\mathfrak m}^{n+2}}$.

De mani\`ere similaire, soit $\tilde \beta\in \mathbb C[[x]]$ une
s\'erie telle que $e^{\beta_0}=\gamma_0$
et $\hbox{Exp}(\tilde \beta;\tilde \beta)\equiv 
E_nx^n\pmod{{\mathfrak m}^{n+1}}$.
Alors la s\'erie 
$$\overline \beta=\tilde
\beta+\frac{1}{\gamma_0}
E_n\left(\sum_{k=1}^{n}\frac{(n\beta_0)^{k-1}}{k!}\right)^{-1}
x^n$$
v\'erifie $\gamma-\hbox{Exp}(\overline \beta;\overline \beta)\equiv 0
\pmod{{\mathfrak m}^{n+2}}$.  
\end{rem}

\begin{rem} Pour $(C,\gamma)\in\mathcal{SI}$, les s\'eries
  $\alpha,\beta\in\mathfrak m$ telles que
$C=\hbox{Exp}(\alpha;\beta),\gamma=x\ \hbox{Exp}(\beta;\beta)$ peuvent
\'egalement se calculer en prenant la s\'erie g\'en\'eratrice de la
premi\`ere colonne des logarithmes matriciels
$$\sum_{n=1}^\infty (-1)^{n+1}\frac{(\rho(C,\gamma)-\hbox{id})^n}{n}
\hbox{ et }
\sum_{n=1}^\infty (-1)^{n+1}\frac{(\rho(\gamma/x,\gamma)-\hbox{id})^n}{n}
$$
des matrices triangulaires inf\'erieures unipotentes $\rho(C,\gamma)$
et $\rho(\gamma/x,\gamma)$.
\end{rem}

{\bf Exemples} La s\'erie formelle 
$\beta\in\mathfrak m$ telle que
$\hbox{Exp}(\beta;\beta)=1+x$ commence par
$$\beta=x-x^2+\frac{3}{2}x^3-\frac{8}{3}x^4+\frac{31}{6}x^5-\frac{157}{15}x^6+\frac{649}{30}x^7+\dots=\sum_{n=0}^\infty
(-1)^nA_n\frac{x^{n+1}}{n!}$$
avec $A_0,A_1,\dots$ donn\'es par
$1,1,3,16,124,1256,15576,226248,3729216,
\dots$, voir la suite A5119 dans \cite{EIS} et la formule (109) dans 
l'article \cite{Labelle} dans lequel la s\'erie formelle 
$x\beta$ est appel\'ee le ''g\'en\'erateur
infinit\'esimal'' de $x\hbox{Exp}(\beta;\beta)\in\mathcal{SD}$
(et qui indique dans la formule (109) \'egalement les premiers
termes de $x\hbox{Exp}(s\beta;s\beta)$).

On peut \'egalement d\'eterminer $\beta(x)$ en utilisant 
l'\'equation diff\'erentielle $xZ'\beta(x)=Z\beta(Z)$ (voir la
preuve de la proposition \ref{propeqdiffbetabeta}) qui se simplifie
en
$$(1+2x)\beta(x)=(1+x)\beta(x+x^2)$$
pour $Z=x(1+x)=x+x^2=x\ \hbox{Exp}(\beta;\beta)$.

Les coefficients de la matrice 
$$\hbox{exp}(u_\beta+d_\beta)=\rho(1+x,x+x^2)=
\left(\begin{array}{rrrrrrr}
1\\1&1\\&2&1\\&1&3&1\\&&3&4&1\\&&1&6&5&1\\
&&&&&&\ddots\end{array}\right)\in {\mathcal I}$$ 
sont donn\'es par
$$\left(\rho(1+x,x+x^2)\right)_{i,j}=
{j+1\choose i-j}=[x^i]\left((1+x)(x+x^2)^j
\right),\ 0\leq i,j\ .$$
La matrice inverse
$$\rho\left(-\frac{1-\sqrt{1+4x}}{2x},-\frac{1-\sqrt{1+4x}}{2}\right)=
\left(\begin{array}{rrrrrrr}
1\\-1&1\\2&-2&1\\
-5&5&-3&1\\
14&-14&9&-4&1\\
-42&42&-28&14&-5&1\\
\vdots&&&&&&\ddots\end{array}\right)$$
avec coefficients non-nuls $(A^{-1})_{i,j}=(-1)^{i+j}
{2i-j\choose i-j}-{2i-j\choose
  i-j-1}=(-1)^{i+j}\frac{j+1}{i+1}{2i-j\choose i-j}$ pour $0\leq j\leq i$
fait intervenir la fonction g\'en\'eratrice
$$\frac{1-\sqrt{1-4x}}{2x}=\sum_{n=0}^\infty \frac{(2n)!}{n!\
  (n+1)!}x^n=1+x+2x^2+5x^3+14x^4+\dots$$
des nombres de Catalan et 
s'obtient \`a un signe pr\`es en lisant les coefficients 
du triangle de Catalan
$$\begin{array}{rrrrrrrrrrrr}
1&\\&1\\1&&1\\&2&&1\\2&&3&&1\\&5&&4&&1\\5&&9&&5&&1\\
&14&&14&&6&&1\\
14&&28&&20&&7&&1\end{array}$$
(avec coefficients $C_{i,j},\ 0\leq i,j$ donn\'es par
${i\choose \frac{i-j}{2}}-{i\choose \frac{i-j}{2}-1}$ si $i\equiv
j\pmod 2$ et $C_{i,j}=0$ sinon) le long de droites affines de pente $1$. 

Pour $\hbox{Exp}(\beta;\beta)=e^x$ on trouve
$$\beta=x-\frac{1}{2}x^2+\frac{5}{12}x^3-\frac{5}{12}x^4
+\frac{107}{240}x^5-\frac{173}{360}x^6+\frac{7577}{15120}x^7+\dots,$$
en accord avec la formule (109)' (qui donne les premiers termes 
de $x\beta$) dans \cite{Labelle}.

\section{Convergence}\label{sectconvergnce}

\begin{thm} \label{thmconv}
Soient $\alpha,\beta\in\mathbb C[[x]]$ deux s\'eries 
formelles d\'efinissant des fonctions holomorphes dans un ouvert
connexe (mais pas n\'ecessairement simplement connexe)
$\mathcal O$ contenant l'origine. Alors il
existe un ouvert $\mathcal V\subset\mathbb C^2$ contenant
$\mathcal O\times \{0\}$ et $\{0\}\times \mathbb C$
tel que $\hbox{Exp}(s\alpha;s\beta)$ est holomorphe pour
$(x,s)\in\mathcal V$.
\end{thm}

\begin{cor} \label{sectconvcor}
Pour $\alpha,\beta\in\mathbb C[[x]]$ deux 
fonctions holomorphes dans un ouvert
connexe
$\mathcal O$ contenant l'origine, le d\'eveloppement en s\'erie en
un point $(\xi,0)\in\mathcal O\times\{0\}$ de la fonction holomorphe 
$\hbox{Exp}(s\alpha;s\beta)\circ (x+\xi)$ 
est donn\'ee par 
$$\sum_{n=0}^\infty G_n\frac{s^n}{n!}$$
o\`u $G_0=1$ et $G_{n+1}=\lbrace\alpha\circ(x+\xi)\rbrace G_n+
(x+\xi)\lbrace\beta\circ(x+\xi)\rbrace G_n'$
o\`u $\alpha\circ(x+\xi)$ et $\beta\circ(x+\xi)$ sont 
les d\'eveloppements en s\'erie  au point 
$\xi\in\mathcal O$ des fonctions holomorphes $\alpha$ et $\beta$.
\end{cor}

Pour $K\in\mathbb C^*$ et $a\in\mathbb C\setminus\{-1\}$ deux
nombres complexes, consid\'erons les s\'eries
formelles $G_0,G_1,G_2,\dots\in \mathbb C[[x]]$
d\'efinies r\'ecursivement par $G_0=1$ et
$G_{n+1}=\left(\left(\frac{K}{1-Kx}\right)^aG_n\right)'$. Le
r\'esultat suivant sera utile dans la preuve du th\'eor\`eme
\ref{thmconv}.

\begin{lem} \label{lemholom} On a
$$\sum_{n=0}^\infty G_n\frac{s^n}{n!}=\left(\frac{1}{1-s(a+1)
\left(\frac{K}{1-Kx}\right)^{a+1}}\right)^{a/(a+1)}$$
et cette s\'erie converge 
pour $(x,s)\in\mathbb C^2$ tel que 
$\vert x\vert<\frac{1}{K}\left(1-K\left((a+1)\vert 
s\vert\right)^{1/(a+1)}\right)$.
\end{lem}

{\bf Preuve} Une r\'ecurrence sur $n$ montre l'\'egalit\'e
$$G_n=\left(\prod_{j=1}^n(j(a+1)-1)\right)\left(\frac{K}{1-Kx}
\right)^{n(a+1)}\ .$$
Le th\'eor\`eme binomial $(1+x)^\alpha=\sum_{n=0}^\infty
  {\alpha\choose n}x^n=\sum_{n=0}^\infty \frac{\alpha(\alpha-1)\cdots
    (\alpha-n+1)}{n!}x^n$ appliqu\'e \`a l'identit\'e
$$\left(\prod_{j=1}^n(j(a+1)-1)\right)=n!(-1-a)^n{-a/(a+1)\choose n}$$
montre
$$\sum_{n=0}^\infty G_n\frac{s^n}{n!}=\left(\frac{1}{1-s(a+1)
\left(\frac{K}{1-Kx}\right)^{a+1}}\right)^{a/(a+1)}$$
avec convergence de la s\'erie pour $(x,s)\in\mathbb C^2$ tel
que $\left\vert s(a+1)\left(\frac{K}{1-Kx}\right)^{a+1}\right\vert<1$.
\hfill$\Box$

{\bf Preuve du th\'eor\`eme \ref{thmconv}} 
Montrons d'abord que la s\'erie 
$\hbox{Exp}(s\alpha;s\beta)$ d\'efinit une fonction 
holomorphe pour $(x,s)$ proche de $(0,0)$.

La proposition \ref{propformExp} montre qu'il suffit pour cela 
de prouver l'analyticit\'e en $(0,0)$ d'une fonction
$\hbox{Exp}(s\tilde\alpha;s\tilde\beta)$ avec 
$\tilde\alpha=\sum_{n=0}^\infty \tilde\alpha_n x^n,
\tilde\beta=\sum_{n=0}^\infty \tilde\beta_n x^n\in\mathbb R[[x]]$ 
des fonctions
holomorphes \`a l'origine telles que $\vert \alpha_n\vert\leq 
\tilde\alpha_n,\ \vert\beta_n\vert\leq \tilde\beta_n$ pour tout $n
\in\mathbb N$.
Comme $\alpha,\beta$ sont holomorphes en $0$, on peut trouver 
$K\geq 1$ tel que $\vert \alpha_n\vert\leq (n+1)K^{n+2},
\vert\beta_n\vert \leq K^{n+1}$
pour tout $n\in\mathbb N$. Il suffit donc de montrer que la fonction 
$\hbox{Exp}\left(s\left(\frac{K}{1-Kx}\right)^2;s\frac{K}{1-Kx}\right)$ 
est analytique en $(0,0)$.
Les s\'eries associ\'ees $g_0,\dots,g_{n+1}=
\left(\frac{K}{1-Kx}\right)^2g_n+x\frac{K}{1-Kx}g_n',\dots$
n'ont que des coefficients positifs et l'in\'egalit\'e $K\geq 1$ 
implique que leurs coefficients sont major\'es par les coefficients
des s\'eries $G_0=1,\dots,G_{n+1}=\left(\frac{K}{1-Kx}G_n\right)'=
\left(\frac{K}{1-Kx}\right)^2G_n+\frac{K}{1-Kx}G_n',\dots$.
Le lemme \ref{lemholom} 
permet donc de minorer le rayon de convergence de
$\hbox{Exp}(s\alpha;s\beta)$ et assure l'analyticit\'e en $(0,0)$ de
la fonction $\hbox{Exp}(s\alpha;s\beta)$.

Pour un point $(\xi,0)$ proche de $(0,0)$, 
la proposition \ref{propformExp} montre par continuation analytique
l'\'egalit\'e
$\hbox{Exp}(s\alpha;s\beta)\circ(x+\xi)=\sum_{n=0}^\infty
  g_{n,\xi}\frac{s^n}{n!}$ 
avec $g_{0,\xi}=1,\dots,g_{n+1,\xi}=\lbrace\alpha\circ(x+\xi)\rbrace
g_{n,\xi}+(x+\xi)\lbrace \beta\circ(x+\xi)\rbrace g_{n,\xi}'$
pour le d\'eveloppement en s\'erie au voisinage d'un point 
$(\xi,0)$ convenable. La minoration du rayon de convergence
de cette s\'erie par le lemme \ref{lemholom} montre que ce 
d\'eveloppement est valable pour tout $\xi$ dans l'ouvert connexe
$\mathcal O$.

Pour montrer l'analyticit\'e de $\hbox{Exp}(s\alpha;s\beta)$
dans un voisinage ouvert de $\mathbb C\times \{0\}$, on utilise
la loi de groupe sur $\mathcal I$ en appliquant it\'erativement 
l'assertion (iv) de la proposition \ref{idfondexp}.
\hfill$\Box$

Le corollaire \ref{sectconvcor} a \'et\'e d\'emontr\'e au cours de 
la preuve du th\'eor\`eme \ref{thmconv}.

\begin{rem} On peut \'egalement donner une preuve du th\'eor\`eme
\ref{thmconv} et du corollaire \ref{sectconvcor}
 en utilisant les majorations de la remarque
\ref{pasdedivpourcoeffdeExp}.
\end{rem}

\section{Une application \`a la combinatoire 
\'enum\'erative}\label{sectintcomb}

Pour $W=\sum_{n=0}^\infty W_n\frac{X^n}{n!}\in\mathbb C[[X]]$
une s\'erie formelle, consid\'erons la s\'erie formelle
$Z=\sum_{n=0}^\infty G_n\frac{s^n}{n!}$
o\`u la suite $G_0,G_1,\dots\in \mathbb C[[X]]$ est d\'efinie
r\'ecursivement par $G_0=1$ et $G_{n+1}=(WG_n)'$. 

Notons $Z_0=\sum_{n=0}^\infty
[X^0]G_n\frac{s^n}{n!}\in\mathbb C[[s]]$ l'\'evaluation en $X=0$
de $Z$. Dans ce chapitre nous d\'ecrivons d'abord
une interpr\'etation combinatoire des coefficients
de $Z_0$ en exprimant $Z_0$ comme fonction de partition 
d'un mod\`ele \`a spins. Ensuite nous montrons qu'une primitive de
$Z_0$ est formellement solution d'une \'equation diff\'erentielle.

\begin{rem} \label{remW} Pour $W\in X\ \mathbb C[[X]]$, la proposition 
\ref{propeqdiffbetabeta} montre qu'on a 
$$Z=\hbox{Exp}\left(sW';s\frac{W}{X}\right)=\frac{d}{dX}
\left(X\ \hbox{Exp}\left(s\frac{W}{X};s\frac{W}{X}\right)\right)$$
(la derni\`ere \'egalit\'e r\'esulte de l'assertion (iii), proposition
\ref{idfondexp}).
\end{rem}

\subsection{Une interpr\'etation combinatoire de $Z_0$}

Pour d\'ecrire l'interpr\'etation combinatoire de $Z_0$, nous
utilisons la s\'erie $W$ pour associer une ``\'energie'' \`a
chaque \'el\'ement, appel\'e \'etat, d'un ensemble fini
correspondant \`a l'espace de phase et repr\'esentant tous les
\'etats possibles d'un syst\`eme ``physique'' discret (qui est
purement math\'ematique dans notre cas). Nous empruntons seulement
un peu de vocabulaire
\`a la physique statistique sans entrer dans les d\'etails.
Plus pr\'ecis\'ement, le coefficient $[s^n]Z_0$
de la fonction $Z_0$ compte le nombre de certaines fonctions
de $\{1,\dots,n\}$ dans $\{1,\dots,n\}$ avec une multiplicit\'e
d\'etermin\'ee par $W$. Ces fonctions
sont en bijection avec des arbres enracin\'es 
\`a $n+1$ sommets \'etiquet\'es de mani\`ere croissante 
\`a partir de la racine par l'ensemble $\{0,\dots,n\}$.

Nous d\'esignons par $\{n\}$ l'ensemble fini $\{1,\dots,n\}$ des
$n$ entiers naturels entre $1$ et $n$ et notons $\{n\}^{\{n\}}$ l'ensemble
des $n^n$ fonctions de $\{n\}$ dans $\{n\}$. 

Consid\'erons une suite $W_0,W_1,W_2,\dots\subset \mathbb C$
de s\'erie g\'en\'eratrice exponentielle 
$W=\sum_{n=0}^\infty W_n\frac{X^n}{n!}$ et associons l'\'energie
$$e_W(f)=\prod_{k=1}^n W_{\sharp(f^{-1}(k))}\in \mathbb C$$
\`a une fonction $f\in\{n\}^{\{n\}}$.

\begin{rem} Un petit calcul, laiss\'e au lecteur, montre l'\'egalit\'e
$$\sum_{f\in\{n\}^{\{n\}}} e_W(f)=[X^0]\frac{d^n}{dx^n}\left(W^n\right)\
.$$
\end{rem}

Appelons une fonction $f\in\{n\}^{\{n\}}$
{\it admissible} si $f(k)\leq k$ pour tout $k\in\{n\}$.
Notons $\mathcal A_n\subset \{n\}^{\{n\}}$ le sous-ensemble de 
toutes les fonctions admissibles. 
Il est facile de voir que $\mathcal A_n$
contient exactement $n!$ \'el\'ements. Plus pr\'ecis\'ement, 
l'application $f\in\mathcal A_n\longmapsto \sum_{k=1}^n (f(k)-1)k!$
est une bijection entre $\mathcal A_n$ et les $n!$ entiers
$0,\dots,(n!-1)$.

\begin{rem} 
La construction suivante donne une bijection classique entre 
l'ensemble $\mathcal A_n$ des fonctions admissibles
et les permutations de $\{n\}$ : 
Pour $f\in\mathcal A_n\subset\{n\}^{\{n\}}$
une fonction admissible, posons $f_1=f$
et d\'efinissons $f_2,\dots,f_n\in  \{n\}^{\{n\}}$ r\'ecursivement
comme suit: $f_{k}(i)=f_{k-1}(i)$ si $i\geq k$ ou si on a 
l'in\'egalit\'e $f_{k-1}(i)<f_{k-1}(k)$ pour $i<k$.
Dans les autres cas, c'est-\`a-dire pour $i<k$ tel que
$f_{k-1}(i)\geq f_{k-1}(k)$, on pose $f_k(i)=f_{k-1}(i)+1$.

On montre facilement que $f_n$ est une permutation de $\{n\}$ et que
la restriction de l'application $f\longmapsto f_n$ \`a l'ensemble
$\mathcal A_n\subset\{n\}^{\{n\}}$ des fonctions admissibles 
est injective.
\end{rem}

L'ensemble $\mathcal A_n$ peut \'egalement \^etre consid\'er\'e
comme l'ensemble des arbres enracin\'es \`a $n+1$ sommets
num\'erot\'es de $0$ \`a $n$ de fa\c con croissante depuis la
racine (on a donc $v>w$ si le sommet num\'erot\'e $v$ est
un successeur du sommet num\'erot\'e $w$).
En effet, associons \`a une fonction admissible $f\in\mathcal A_n$
l'arbre $T_f$ dont la racine porte le num\'ero $0$. 
Le pr\'ed\'ecesseur imm\'ediat d'un sommet $j>0$ 
est le sommet $f(j)-1$. Il est facile de voir que
l'application $f\in\mathcal A_n\longmapsto T_f$ est bijective.
De plus, elle pr\'eserve les ``degr\'es'': Soit $\mathcal D_f(k)
=\{i\in \{n\}\ \vert\  \sharp(f^{-1}(i))=k\}
\subset \{n\}$ le sous-ensemble des 
\'el\'ements dans $\{n\}$ ayant $k$ pr\'eimages sous l'application 
$f\in\mathcal A_n$. Alors un sommet
num\'erot\'e $i<n$ de $T_f$ est de degr\'e $k\geq 1$ si et seulement si
$i+1\in\mathcal D_f(k)$ o\`u le degr\'e d'un sommet est 
d\'efini comme le nombre de ses successeurs imm\'ediats. 
En particulier, le nombre de sommets de degr\'e $k\geq 1$ de
$T_f$ est donn\'e par le nombre d'\'el\'ements 
$\sharp(\mathcal D_f(k))$ du
sous-ensemble $\mathcal D_f(k)\subset \{n\}$. Le nombre de feuilles
(sommets de degr\'e $0$) de $T_f$ est donn\'e par 
$1+\sharp(\mathcal D_f(0))$ o\`u $\mathcal D_f(0)$ est l'ensemble
compl\'ementaire $\{n\}\setminus f(\{n\})$ 
de l'image $f(\{n\})\subset \{n\}$.

\begin{rem} Le chapitre 1.3 dans \cite{St1} d\'ecrit une
bijection entre les permutations de $\{n\}$ et les arbres enracin\'es
avec $(n+1)$ sommets num\'erot\'es de $0$ \`a $n$ de mani\`ere
croissante.
\end{rem} 

Consid\'erons maintenant la fonction $Z_0=\sum_{n=0}^\infty \lbrace
[X^0]G_n\rbrace\frac{s^n}{n!}$ d\'efinie \`a partir 
de $G_0=1,\dots,G_{n+1}=(WG_n)',\dots$ pour une s\'erie 
g\'en\'eratrice exponentielle
$W=\sum_{n=0}^\infty W_n\frac{X^n}{n!}$ comme au d\'ebut du chapitre.

\begin{prop} \label{propidenum} On a l'\'egalit\'e
$$Z_0=1+\sum_{n=1}^\infty \frac{s^n}{n!}\sum_{f\in \mathcal
  A_n}e_W(f)\ .$$
\end{prop}

{\bf Exemples} (i) Un petit calcul facile montre 
qu'on obtient $Z_0=\frac{1}{1-s}$ pour $W=e^X$. Cet exemple 
trivial est \'equivalent \`a l'identit\'e $\sharp(\mathcal A_n)=n!$.

\ \ (ii) Pour $W=\left(\frac{K}{1-K(X+\xi)}\right)^a$ avec 
$a\in\{1,2,\dots\},K\in\mathbb C^*$ et $\xi\in \mathbb
C\setminus\{1/K\}$, le lemme \ref{lemholom} donne
$$Z_0=\left(1-s(a+1)\left(\frac{K}{1-K\xi}\right)^{a+1}\right)^{-a/(a+1)}
\ .$$

\ \ (iii) Le choix de
$$W=\frac{(X+i)^2}{1-(X+i)^2}=-1+\sum_{n=0}^\infty(-1)^n\frac{X^{4n}}{4^{n+1}}\left(2+2iX-X^2\right)$$
donne
$$\begin{array}{l}\displaystyle
Z=\frac{d}{dX}\left(\left(-\sqrt{1-2s(X+i)+(s^2-2)(X+i)^2
-2s(X+i)^3+(X+i)^4}\right.\right.\\
\displaystyle \left.+1-s(X+i)+(X+i)^2\Big)\Big/(2(X+i))\right)
\end{array}$$
(voir remarque \ref{remW} et l'exemple (i) du chapitre
\ref{expbetabeta}). On obtient donc
$$Z_0=1+\frac{is}{\sqrt{4-s^2}}=1+\frac{i}{2}\sum_{n=0}^\infty{2n\choose
  n}\frac{s^{2n+1}}{16^n}\ .$$

Notons $\mathcal S(2n,k)$ l'ensemble fini form\'e par tous les 
arbres (non-planaires)
enracin\'es de degr\'e maximal $\leq 2$ (chaque sommet poss\`ede donc
au plus deux descendants directs) avec $2n$ sommets num\'erot\'es de
mani\`ere croissante dans $\{0,\dots,2n-1\}$ et $k$ feuilles
(sommets terminaux).

Introduisons \'egalement l'ensemble fini $\mathcal E(2n+1,k)$ des
arbres (non-planaires) enracin\'es n'ayant que des sommets de
degr\'es pairs (le nombre de descendants directs de chaque sommet
est pair) avec $2n+1$ sommets num\'erot\'es de mani\`ere croissante
dans $\{0,\dots,2n\}$ et $k$ sommets int\'erieurs (de degr\'es $\geq
2$). 

Nous allons donner plus loin une nouvelle preuve du r\'esultat suivant
qui se trouve dans \cite{KPP} : 
 
\begin{thm}\label{idcomben} Pour tout $n,k\in\mathbb N$, les ensembles
finies $\mathcal S(2n,k)$ et $\mathcal E(2n+1,k)$ contiennent le
m\^eme nombre d'\'el\'ements.
\end{thm}

{\bf Preuve de la proposition \ref{propidenum}} 
Notons $W(X_i)=\sum_{n=0}^\infty W_n\frac{X_i^n}{n!}
\in\mathbb C[[X_i]]$ la s\'erie formelle
$W$ en une variable $X_i$ et posons $\tilde G_0=1,\dots,\tilde G_{n+1}=
\left(\sum_{j=1}^{n+1} t_j\frac{\partial}{\partial X_j}\right)
\left(W(X_{n+1})\tilde G_n\right),\dots$. On obtient alors $Z$ 
(respectivement $Z_0$)
en \'evaluant  la s\'erie formelle 
$$\sum_{n=0}^\infty \tilde G_n\frac{s^n}{n!}\in\mathbb C[[s,X_1,X_2,\dots,
t_1,t_2,\dots]]$$
en $t_j=1$ et $X_j=X$ (respectivement en $t_j=1$ et $X_j=0$)
pour tout $j\geq 1$. Les contributions \`a $[X_1^0X_2^0\cdots
X_{n}^0]\tilde G_n$ sont de la forme
$$\sum_{i_1,\dots,i_n\in\{n\}}\alpha_{i_1\dots i_n}\prod_{j=1}^n
(t_j^{i_j}W_{i_j})$$
o\`u $\alpha_{i_1\dots i_n}$ est le nombre de fonctions admissibles
$f\in \{n\}^{\{n\}}$ telles que $\sharp(f^{-1}(j))=i_j$. En effet,
l'exposant de $t_j$ indique l'ordre de la d\'erivation par
rapport \`a la variable $X_j$ et une fonction $f:\{1,2,\dots\}
\longrightarrow\{1,2,\dots\}$ d\'etermine
pour chaque entier $j\leq n$ le
sommand $t_{f(j)}\frac{\partial}{\partial X_{f(j)}}$ 
de l'op\'erateur diff\'erentiel $\sum_{j=1}^\infty
t_j\frac{\partial}{\partial X_j}$ \`a appliquer \`a 
$(W(X_j)\tilde G_{j-1})$ pour obtenir une contribution 
\`a $\tilde G_j$.
Comme seules les variables $X_1,\dots,X_j$ apparaissent dans
$W(X_j)\tilde G_{j-1}$, on peut se restreindre aux fonctions admissibles.
\hfill$\Box$

\subsection{Une \'equation diff\'erentielle pour $Z_0$}

Consid\'erons comme pr\'ec\'edemment la s\'erie formelle
$Z_0\in\mathbb C[[s]]$ obtenue par l'\'evaluation en $X=0$
de $Z=\sum_{n=0}^\infty G_n\frac{s^n}{n!}$ d\'efini \`a partir
de la suite r\'ecurrente $G_0,\dots,G_{n+1}=(WG_n)',\dots\subset\mathbb 
C[[X]]$ pour $W\in\mathbb C[[X]]$ une s\'erie formelle donn\'ee.

\begin{thm} \label{propdiffeqpouru}
Soit $Z_0=\sum_{n=0}^\infty \zeta_n s^n
\in\mathbb C[[s]]$ une s\'erie formelle associ\'ee
\`a $W=\sum_{n=0}^\infty W_n\frac{X^n}{n!}$ 
comme ci-dessus. Alors la primitive $U=\sum_{n=0}^\infty U_ns^n=
\sum_{n=0}^\infty \zeta_n
\frac{s^{n+1}}{n+1}\in s\mathbb C[[s]]$ est une solution 
formelle de l'\'equation diff\'erentielle
$$W_0 U'=W\circ(W_0 U)$$
avec la condition initiale $U(0)=0$.
\end{thm}

\begin{rem} (i) Si $W_0=0$, l'\'equation diff\'erentielle
du th\'eor\`eme \ref {propdiffeqpouru} d\'eg\'en\`ere. 
Un calcul facile montre cependant qu'on a dans ce cas $Z_0=e^{W_1s}$.

\ \ (ii) Le th\'eor\`eme
  \ref{propdiffeqpouru} permet un calcul 
r\'ecursif des coefficients de $U$ par ``bootstrapping'' : 
\`A partir d'une solution approch\'ee \`a l'ordre $n$ de $U$,
l'\'equation diff\'erentielle donne une solution 
approch\'ee \`a l'ordre $n$ de $U'$ ce qui d\'etermine 
$U$ \`a l'ordre $n+1$.

\ \ (iii) En int\'egrant l'\'equation diff\'erentielle
du th\'eor\`eme \ref {propdiffeqpouru} on obtient 
la fonction r\'eciproque 
$$s=W_0\int_0^U\frac{dv}{W\circ(W_0v)}=\int_0^{W_0U}\frac{dv}{W(v)}
\in U\mathbb C[[U]]$$
de $U$.

\ \ (iv) On peut se d\'ebarasser du terme constant $W_0\not=1$
En consid\'erant $\tilde Z_0$ associ\'e \`a $\tilde W=\frac{1}{W_0}
W=1+\dots$ et en remarquant qu'on a $Z_0=\tilde Z_0\circ(W_0 s)$
pour $Z_0$ associ\'e \`a $W$.
\end{rem}

{\bf Preuve du th\'eor\`eme \ref{propdiffeqpouru}}
Consid\'erons d'abord une s\'erie formelle 
$W=\sum_{n=0}^\infty W_n\frac{X^n}{n!}$ avec $W_0\not=0$ 
qui d\'efinit une fonction holomorphe dans un ouvert connexe de 
$\mathbb C$ contenant $0$ et une racine $\xi$ de $W$.
La fonction $w(x)=W(x+\zeta)$ est alors holomorphe
dans un ouvert connexe $\mathcal O\subset \mathbb C$ 
contenant $0$ et $-\xi$. Comme $w$ admet 
l'origine $0\in\mathbb C$ comme racine, 
la s\'erie $\frac{w}{x}\in \mathbb C[[x]]$ d\'efinit 
\'egalement une fonction holomorphe
dans $\mathcal O$. Le th\'eor\`eme \ref{thmconv}
montre que les fonctions $\hbox{Exp}(sw';s\frac{w}{x})$ 
et $\hbox{Exp}(s\frac{w}{x};s\frac{w}{x})$
sont holomorphes en $x$ et $s$ pour $(x,s)$ appartenant \`a 
un ouvert contenant $\mathcal O\times \{0\}\subset \mathbb C^2$.

Consid\'erons maintenant les fonctions holomorphes
$$y_s=\hbox{Exp}\left(sw';s\frac{w}{x}\right)\vert_{x=-\xi},\
z_s=\hbox{Exp}\left(s\frac{w}{x};s\frac{w}{x}\right)\vert_{x=-\xi}$$
de $s$, obtenues en \'evaluant $\hbox{Exp}(sw';s\frac{w}{x})$ 
et $\hbox{Exp}(s\frac{w}{x};s\frac{w}{x})$ en $x=-\xi$. 
Comme $w'=\left(x\frac{w}{x}\right)'$, 
l'assertion (iii) de la proposition
\ref{idfondexp} montre l'identit\'e 
$$y_s=\left(\frac{d}{dx}(x\hbox{Exp}
(s\frac{w}{x};s\frac{w}{x}))\right)\vert_{x=-\xi}\ .$$
L'\'equation diff\'erentielle de l'assertion (ii), proposition
\ref{diffeqpourExp}, devient donc
$$\begin{array}{lcl}
y_s'&=&y_s\ w'\circ(-\xi z_s)\\
z_s'&=&-\frac{w(-\xi)}{\xi}y_s\ .\end{array}
$$
En \'eliminant $y_s=-\frac{\xi}{w(-\xi)}z_s'=-\frac{\xi}{W_0}z_s'$,
nous obtenons
$$z_s''=z_s'\ w'\circ(-\xi z_s)=z_s'\ W'\circ(\xi-\xi z_s)$$
avec les conditions initiales $z_0=1$ et $z_0'=-\frac{W_0}{\xi}y_0=
-\frac{W_0}{\xi}$. La primitive $U=\frac{1}{W_0}(\xi-\xi z_s)$
de $y_s$ satisfait alors l'\'equation diff\'erentielle 
$U''=U'\ \{W'\circ(W_0U)\}=\frac{1}{W_0}\big(W\circ(W_0U)\big)'$ 
du th\'eor\`eme pour les conditions initiales
$U(0)=0$ et $U'(0)=1$. Une int\'egration donne alors
$U'=\frac{1}{W_0}W\circ(W_0U)$.

Pour d\'emontrer le cas g\'en\'eral, il suffit de prouver la
congruence $U'=W\circ(W_0U)\pmod{s^N}$ pour $N\in\mathbb N$
un entier naturel arbitraire. Fixons un tel entier $N\in\mathbb
N$ et consid\'erons le polyn\^ome $\tilde W=X^{N+1}+\sum_{n=0}^{N}
W_n\frac{X^n}{n!}$. Un petit calcul montre que la s\'erie
$\tilde Z_0$ associ\'ee \`a $\tilde W$ coincide avec la s\'erie 
$Z_0$ associ\'ee \`a $W$ jusqu'\`a l'ordre $N$.
Comme $\tilde W$ est holomorphe dans $\mathbb C$ 
et poss\`ede au moins une racine, la preuve pr\'ec\'edente 
marche pour $\tilde U=\int \tilde Z_0$. Le r\'esultat 
suit maintenant des observations $W\equiv \tilde W\pmod{X^{N}}$
et $U=\int Z_0\equiv  \tilde U=\int \tilde Z_0\pmod{s^{N+1}}$
avec $U,\tilde U\in s\mathbb C[[s]]$.
\hfill$\Box$ 

\subsection{Exemples et preuve du th\'eor\`eme \ref{idcomben}}

\subsubsection{Exemple}\label{explerationnel} Pour $W=\frac{1}{P}$
avec $P\in \mathbb C[X]$ un polyn\^ome tel que $P(0)=1/W_0\not=0$, on a
$$\int_0^{W_0U} P(v)dv =s\ .$$
En particulier, $U$ et $Z_0=U'$ sont alg\'ebriques. Cette 
construction est \'evidemment reli\'ee au chapitre
\ref{expbetabeta}. On peut ainsi \'egalement consid\'erer
une fraction rationelle $\beta\in\mathfrak m$ satisfaisant la
condition du corollaire \ref{expalg} (qui est \'equivalente
au fait que la fraction rationelle $\frac{1}{x\beta}$ 
admet une primitve rationelle). Le choix d'un nombre complexe
$\xi\in\mathbb C$ tel que $\beta(\xi)\in\mathbb C^*$ 
donne alors une fonction alg\'ebrique $Z_0$ en consid\'erant 
la $W(X)=\frac{\beta(X-\xi)}{X-\xi}$.  

Similairement, pour $H=\sum_{n=0}^\infty H_nX^n\in\mathbb C[[X]]$
une s\'erie formelle telle que $H_0=H(0),H_1=H'(0)\not=0$,
on peut consid\'erer $W=\sum_{n=0}^\infty W_n\frac{X^n}{n!}=
\frac{H}{H'}$ avec $W(0)=W_0\not=0$ et on obtient $H(W_0U)=H(0)e^s$. 

\subsubsection{Exemple}\label{explearbre}
Pour le polyn\^ome $W=1+tX+\frac{X^2}{2}$ on trouve 
l'\'equation diff\'erentielle
$$U'=1+tU+\frac{U^2}{2}$$
avec la condition initiale $U(0)=0$ admettant la solution
$$U=\hbox{tan}\left(\frac{s}{2}\sqrt{2-t^2}+
\hbox{arctan}\left(\frac{t}{\sqrt{2-t^2}}\right)
\right)\sqrt{2-t^2}-t$$
avec d\'eriv\'ee
$$Z_0=\frac{2-t^2}{2}\left(1+\hbox{tan}\left(\frac{s}{2}
\sqrt{2-t^2}+\hbox{arctan}\left(\frac{t}{\sqrt{2-t^2}}\right)
\right)^2\right)\ .$$

La sp\'ecialisation
$t=1$ donne
$$Z_0=\frac{1+(\hbox{tan}(s/2+\pi/4))^2}{2}=
\frac{1+\hbox{sin}(s)}{\hbox{cos}(s)^2}=1+s+2\frac{s^2}{2!}+
5\frac{s^3}{3!}+16\frac{s^4}{4!}+\dots$$
et correspond \`a la s\'erie g\'en\'eratrice exponentielle
de la suite des nombres $$1,1,2,5,16,61,\dots,$$
aussi appel\'es nombres d'Euler, voir la suite A111 dans \cite{EIS}
ou le chapitre 3.16 dans \cite{St1}.

La sp\'ecialisation $t=0$ donne
$$Z_0=1+(\hbox{tan}(s/\sqrt
2))^2=-2\left(\hbox{log}\circ\hbox{cos}\left(\frac{s}{\sqrt 2}\right)
\right)''=1+\frac{s^2}{2!}+
4\frac{s^4}{4!}+34\frac{s^6}{6!}+\dots$$
et correspond \`a la s\'erie g\'en\'eratrice exponentielle 
paire de la suite enti\`ere
$$1,1,4,34,496,\dots$$
reli\'ee aux nombres de Bernoulli, voir la suite A2105 dans \cite{EIS}.

Pour $W=1+tX+\frac{X^2}{2}$ avec $t$ non \'evalu\'e, le coefficient 
$A_n=[s^n]Z_0\in \mathbb N[t]$ est un polyn\^ome de degr\'e $n$
en $t$, appel\'e le $n-$i\`eme polyn\^ome d'Andr\'e, voir 
la suite A94503 dans \cite{EIS}. Le coefficient $[t^k]A_n$
donne le nombre de fonctions admissibles dans $\mathcal A_n$
telles que chaque \'el\'ement $j\in\{n\}$ a au plus $2$
pr\'eimages et exactement $k$ \'el\'ements
de $\{n\}$ ont une pr\'eimage unique. 
Les polyn\^omes d'Andr\'e sont faciles \`a calculer
en utilisant le r\'esultat suivant.

\begin{lem} La suite $A_0=1,A_1=t,A_2=1+t^2,A_3=4t+t^3,\dots$ 
des polyn\^omes d'Andr\'e satisfait la relation de r\'ecurrence
$$A_{n+1}=(1-t^2/2)A_n'+t\frac{n+2}{2}A_n\ .$$
\end{lem}

{\bf Preuve} Soit $f\in \mathcal A_{n+1}$ une fonction admissible
donnant une
contribution de $1$ \`a $[t^k]A_{n+1}$. Consid\'erons la restriction
$g$ de $f$ \`a $\{n\}$. Comme $f\in\{n+1\}^{\{n+1\}}$ est admissible,
$g\in \{n\}^{\{n\}}$ l'est aussi et fournit une contribution de 
$1$ \`a $[t^{k+1}]A_n$ si la pr\'eimage de $f(n+1)$ est double
ou une contribution de $1$ \`a $[t^{k-1}]A_n$ si la pr\'eimage
de $f(n+1)$ est simple.
Dans le premier cas, il y a $(k+1)$ possibilit\'es pour choisir
$f(n+1)$. Dans le deuxi\`eme cas, $g$ a une pr\'eimage double pour
$a$ \'el\'ements dans $\{n\}$ avec $n=2a+(k-1)$ ce qui donne 
$a=\frac{n+1-k}{2}$ et le 
nombre de choix possibles pour $f(n+1)$ est donn\'e par
$n+1-(k-1)-a=\frac{n-k+3}{2}=\frac{n+2-(k-1)}{2}$.
On a donc la r\'ecurrence $A_{n+1}=\left(1-\frac{t^2}{2}\right)A_n'+
t\frac{n+2}{2}A_n$.\hfill $\Box$


\subsubsection{Exemple}\label{explepair}
En posant $W=1-t+t\ \hbox{cosh}(X)=1+t\sum_{n=0}^\infty
\frac{X^{2n}}{(2n)!}$ on obtient l'\'equation diff\'erentielle
$$U'=1-t+t\ \hbox{cos}(U)$$
avec la condition initiale $U(0)=1$
satisfaite par la fonction 
$$\begin{array}{l}\displaystyle 
U=\frac{1}{2}\left(\hbox{log}\left(
1+\left(\hbox{tan}\left(\frac{s}{2}\sqrt{2t-1}+\hbox{arctan}\left(\frac{1}{\sqrt{2t-1}}\right)\right)\right)^2\right)\right.\\
\displaystyle \left. \qquad -\hbox{log}\left(
1+\left(\hbox{tan}\left(-\frac{s}{2}\sqrt{2t-1}+\hbox{arctan}\left(\frac{1}{\sqrt{2t-1}}\right)\right)\right)^2\right)\right)\
.\end{array}$$
La d\'eriv\'ee $Z_0=U'$ est donc donn\'ee par 
$$\begin{array}{l}\displaystyle 
Z_0=\frac{\sqrt{2t-1}}{2}\left(\hbox{tan}
\left(\frac{s}{2}\sqrt{2t-1}+\hbox{arctan}\left(\frac{1}{\sqrt{2t-1}}
\right)\right)\right.\\
\displaystyle \qquad \left. +\hbox{tan}
\left(-\frac{s}{2}\sqrt{2t-1}+\hbox{arctan}\left(\frac{1}{\sqrt{2t-1}}
\right)\right)\right)\ .\end{array}$$

{\bf Preuve du th\'eor\`eme \ref{idcomben}} Il suffit de v\'erifier 
que la substitution $s\longmapsto st,t\longmapsto \frac{1}{t^2}$
dans la d\'eriv\'ee $Z_0'$ de l'exemple \ref{explepair}
coincide (\`a un facteur $t^2$ pr\`es) avec la partie
impaire en $s$ de la fonction $Z_0$ apparaissant dans
l'exemple \ref{explearbre}. Nous laissons les d\'etails calculatoires
au lecteur. \hfill$\Box$

\section{Matrices de Riordan}\label{sectriordan}

L. Shapiro a introduit le groupe matriciel $\rho(\mathcal I)$
sous le nom de ``groupe de Riordan'', voir par exemple
\cite{SGWW}. Encore auparavant, certains
\'el\'ements de ce groupe ont \'et\'e \'etudi\'e par D.G. Rogers sous
le nom de ``renewal arrays'', voir \cite{R}. Ce chapitre est un bref
r\'esum\'e des travaux entrepris par diff\'erents auteurs
dans le but d'\'elucider certaines structures combinatoires \`a
l'aide de matrices de Riordan.

Commen\c cons par citer le r\'esultat suivant, d\^u \`a D.G. Rogers,
utile pour caract\'eriser les matrices dans $\rho(\mathcal I)$.

\begin{thm} \label{thmrogers}
(D.G. Rogers, voir \cite{R}?) Soit $M$ une matrice
  triangulaire inf\'erieure infinie de coefficients $M_{i,j},0\leq
  i,j$. Alors $M\in \rho(\mathcal I)$ si et seulement si il existe
une suite $a_0,a_1,\dots$ de nombres complexes avec $a_0\not= 0$
telle que $M_{n+1,k+1}=\sum_{j=0}^{n-k}\alpha_jM_{n,k+j}$.
\end{thm} 
 
\begin{rem} Pour $\rho(B,\beta)\in \rho(\mathcal I)$,
on a $\sum_{n=0}^\infty a_nx^n=\frac{x}{\beta^{\langle -1\rangle}}$ pour
la s\'erie g\'en\'eratrice $A=\sum_{n=0}^\infty a_nx^n$ 
de la suite $A=(a_0,a_1,\dots)$
(appel\'ee $A-$suite de la matrice de Riordan associ\'ee
\`a $(B,\frac{\beta}{x})$ dans la litt\'erature sur les matrices 
de Riordan) intervenant dans le th\'eor\`eme \ref{thmrogers}.

Mentionnons \'egalement l'existence d'une 
g\'en\'eralisation de ce th\'eor\`eme, voir \cite{MRSV}. 
\end{rem}

{\bf Preuve du th\'eor\`eme \ref{thmrogers}} 
L'identit\'e facile 
$$\left(B\left(\frac{\beta}{x}\right)^n,\beta\right)
\left(\frac{x}{\beta^{\langle -1\rangle}},x\right)=
\left(B\left(\frac{\beta}{x}\right)^{n+1},\alpha\right)$$
montre le r\'esultat pour un \'el\'ement
$\rho(B,\beta)\in\rho(\mathcal I)$
en consid\'erant la suite 
$a_0,a_1,\dots $ des coefficients de $\frac{x}{\beta^{\langle
    -1\rangle}}=\sum_{n=0}^\infty a_nx^n$.

D'autre part, soit $B=\sum_{n=0}^\infty B_nx^n$ la s\'erie 
g\'eneratrice de la premi\`ere colonne d'une matrice $M$
satisfaisant l'hypoth\`ese
du th\'eor\`eme \ref{thmrogers} pour une suite $a_0,a_1,\dots$. 
Consid\'erons la matrice $\tilde M=\rho(B,\beta)$ pour
$\beta=\left(\frac{x}{\sum_{n=0}^\infty a_nx^n}\right)^{\langle
  -1\rangle}$. La premi\`ere colonne des matrices $M$ et $\tilde M$
  coincide alors et les autres colonnes de $M$ et $\tilde M$ sont
d\'efinies r\'ecursivement de la m\^eme mani\`ere. On a donc 
$M=\tilde M$.\hfill$\Box$

La formule de r\'ecursion pour les coefficients donn\'ee par le
th\'eor\`eme \ref{thmrogers} peut s'interpr\'eter en termes 
de chemins sur le r\'eseau $\mathbb Z^2$ comme suit : Consid\'erons 
les chemins reliant l'origine $(0,0)$ au point $(n,k)$ sans quitter
le c\^one $\{k\geq 0\}\cap \{n\geq k\}$ et utilisant des pas 
color\'es par $a_j$ couleurs de la forme $(1,1-j)$ pour
$j=0,1,2,\dots$. Le nombre $m_{n,k}$ de tels chemins 
satisfait par construction l'\'equation 
$$m_{n+1,k+1}=\sum_{j=0}^{n-k} a_jm_{n,k+j}\ .$$
Le th\'eor\`eme \ref{thmrogers} montre donc pour $a_0\not= 0$ 
l'appartenance \`a $\rho(\mathcal I)$ de la matrice triangulaire
inf\'erieure avec coefficients $m_{i,j},\ 0\leq i,j$.
Plus pr\'ecis\'ement, cette matrice est donn\'ee par
$\rho\left(\frac{\beta}{x},\beta\right)$ avec
$\beta=\big(x/\sum_{n=0}^\infty a_n x^n\big)^{\langle -1\rangle}$.

{\bf Exemples} Pour $a_0a_2\not=0$ et $a_3=a_4=\dots=0$, on trouve
$$\beta=\frac{(1-a_1x)-\sqrt{(1-a_1x)^2-4a_0a_2}}{2a_2x}\ .$$
Des choix convenables pour $a_0,a_1,a_2$ permettent de retrouver
des triangles apparaissant dans la lit\'erature de la combinatoire 
\'enum\'erative.

Mentionnons \'egalement l'\'el\'ement $\rho(1,e^x-1)\in\mathcal I$
qui est reli\'e aux nombres de Stirling de deuxi\`eme esp\`ece,
voir par exemple le chapitre 5.8 de \cite{Spr} pour les d\'etails.

Pour terminer ce chapitre nous relions bri\`evement ce qui pr\'ec\`ede 
\`a certains d\'eterminants.

Associons \`a une suite $a_0,a_1,a_2,\dots$ la suite
$d_0=1,d_1=a_0,\dots\subset \mathbb C[x]$ des 
d\'eterminants $d_n=\det(T_n)$
o\`u $T_n$ est la matrice de Toeplitz
$$T_n=\left(\begin{array}{ccccccc}
a_0&-x&0&\dots&0\\
a_1&a_0&-x&&\vdots\\
a_2&a_1&a_0&&\\
\vdots&&&\ddots&-x\\
a_{n-1}&a_{n-2}&\dots&a_1&a_0\end{array}\right)$$
avec coefficients $T_{i,j}=a_{i-j}$ si $0\leq j\leq i<n$, 
$T_{i,i+1}=-x$ et $T_{i,j}=0$ sinon.
Un calcul facile donne la r\'ecurrence
$$d_{n+1}=\sum_{j=0}^n a_jd_{n-j}x^j\ .$$
La suite $d_0,d_1,\dots$ est d'ailleurs \'egalement la suite
des d\'eterminants des matrices $\tilde T_n$ pour 
$$\tilde T_n=\left(\begin{array}{ccccccc}
a_0&-1&0&\dots&0\\
a_1x&a_0&-1&&\vdots\\
a_2x^2&a_1x&a_0&&\\
\vdots&&&\ddots&-1\\
a_{n-1}x^{n-1}&a_{n-2}x^{n-2}&\dots&a_1x&a_0\end{array}\right)\ .$$

Si $a_0\not=0$, le Th\'eor\`eme \ref{thmrogers} montre que la matrices 
triangulaire infinie avec coefficients $M_{i,j}=[x^{i-j}]d_{i+1},\ 0\leq
i,j$ est de la forme $\rho\left(\frac{\beta}{x},\beta\right)\in
\rho(\mathcal I)$ avec $\beta=\left(\frac{x}{\sum_{n=0}^\infty a_nx^n}
\right)^{\langle -1\rangle}$.

\section{L'alg\`ebre ${\mathcal P}$  des matrices triangulaires
infinies polynomiales}\label{sectmatrpolyn}

Ce chapitre d\'ecrit une alg\`ebre contenant quelques sous-alg\`ebres
et groupes int\'eressants. En particulier, elle contient le groupe
d'interpolation $\mathcal I$ consid\'er\'e dans les chapitres 
pr\'ec\'edents.

{\bf D\'efinition} On dira qu'une matrice triangulaire inf\'erieure
infinie $M$ est {\it polynomiale} s'il existe une suite 
de polyn\^omes $p_0,p_1,\dots\in{\mathbb C}[u]$ telle que
$M_{i,j}=p_{i-j}(j)$ pour tout coefficient
non-nul $M_{i,j},0\leq j\leq i$ de $M$. 

D\'esignons par ${\mathcal P}$ l'espace vectoriel 
des matrices triangulaires
inf\'erieures infinies qui sont polynomiales.
Notons $\mathcal{SG}\subset {\mathcal P}$ le sous-ensemble
des matrices unipotentes (avec coefficients $1$ sur la
diagonale) dans ${\mathcal P}$ et notons $\mathfrak{sg}
\subset {\mathcal P}$ le sous-espace vectoriel des matrices
triangulaires inf\'erieures strictes de ${\mathcal P}$.

\'Etant donn\'ee une matrice M$\in\mathcal P$ avec coefficients
$M_{i,j}=p_{i-j}(j)$, notons $\tau_\lambda(M)\in \mathcal P$ la
matrice avec coefficients
$\left(\tau_\lambda(M)\right)_{i,j}=p_{i-j}(j+\lambda)$.
 
\begin{thm} \label{thmalg} (i) L'ensemble ${\mathcal P}$ est 
une alg\`ebre. L'ensemble 
$\mathcal{SG}$ est un groupe de Lie de dimension infinie 
d'alg\`ebre de Lie $\mathfrak{sg}$.

\ \ (ii) Pour tout $\lambda\in\mathbb C$, 
l'application lin\'eaire $\tau_\lambda:\mathcal
P\longrightarrow \mathcal P$ d\'efinie par $M\longmapsto
\tau_\lambda(M)$
est un automorphisme de l'alg\`ebre $\mathcal P$ et de 
l'alg\`ebre de Lie $\mathfrak{sg}$. 
\end{thm}

\begin{rem} (i) La preuve du th\'eor\`eme \ref{thmalg} montre que l'espace
  vectoriel $\mathfrak g=\mathcal P$ est \'egalement une alg\`ebre
de Lie. Le groupe de Lie associ\'e n'est cependant pas dans $\mathcal
P$ car il contient \'egalement les matrices diagonales avec coefficients 
diagonaux $e^{p_0(0)},e^{p_0(1)},e^{p_0(2)},\dots$ pour $p_0\in\mathbb
C[x]$ un polyn\^ome. On pourrait cependant l\'eg\`erement agrandir 
l'alg\`ebre de Lie $\mathfrak{sg}$ en remarquant que les
\'el\'ements inversibles de l'alg\`ebre $\mathcal P$ sont 
exactement les matrices de $\mathcal P$ ayant une diagonale constante
non-nulle. L'alg\`ebre de Lie associ\'ee s'obtient en rajoutant 
les multiples de l'identit\'e \`a l'espace vectoriel $\mathfrak{sg}$.

\ \ (ii) L'automorphisme $\tau_\lambda$ est la g\'en\'eralisation \`a
$\mathcal P$ de
l'automorphisme de groupe $\varphi_{0,\lambda,0}$ de $\mathcal{SI}$
consid\'er\'e dans la proposition \ref{propvarphi}. En particulier,
$\tau_1$ agit en effa\c cant la premi\`ere ligne et
colonne d'une matrice $M\in\mathcal P$.

\end{rem}

L'ingr\'edient crucial pour prouver le th\'eor\`eme \ref{thmalg}
est le r\'esultat suivant. 
 
\begin{lem} \label{lemalg}
(i) Soient $a=(a_0,a_1,\dots,),b=(b_0,b_1,\dots)\in
{\mathbb C}[u]^{\mathbb N}$ deux suites de polyn\^omes. Alors il
existe une suite de polyn\^omes $c=(c_0,c_1,\dots)$ telle que
$$M(a)M(b)=M(c)\ .$$
De plus, on a $\hbox{deg}(c_k)\leq \hbox{sup}_{0\leq h\leq k}
(\hbox{deg}(a_h)+\hbox{deg}(b_{k-h}))$ pour tout $k\in{\mathbb N}$.

\ \ (ii) Soient $a=(a_0,a_1,\dots,),b=(b_0,b_1,\dots)\in
{\mathbb C}[u]^{\mathbb N}$ deux suites de polyn\^omes. Alors les
degr\'es $\hbox{deg}(\tilde c_k)$ de
la suite $\tilde c=(\tilde c_0,\tilde c_1,\dots)$ associ\'ee
au commutateur 
$$M({\tilde c})=[M(a),M(b)]=M(a)M(b)-M(b)M(a)$$
v\'erifient les in\'egalit\'es 
$\hbox{deg}(\tilde c_k)<\hbox{sup}_{0\leq h\leq k}
(\hbox{deg}(a_h)+\hbox{deg}(b_{k-h}))$ pour tout $k\in{\mathbb N}$.
\end{lem}

{\bf Preuve du lemme \ref{lemalg}} Comme le coefficient de
$(M(a)M(b))_{l+k,l}$ ne d\'epend que de $(M(a))_{i,j},(M(b))_{i,j}$
avec $i-j\leq k$, on peut supposer que les
polyn\^omes $a_0,a_1,\dots,b_0,b_1,\dots$ sont tous de degr\'e au plus
$$\hbox{max}(\hbox{deg}(a_0),\dots,\hbox{deg}(a_k),
\hbox{deg}(b_0),\dots,\hbox{deg}(b_k))\ .$$
Par lin\'earit\'e, il suffit de consid\'erer deux
matrices $A,B$ avec coefficients
$A_{i,j}=j^s \alpha_{i-j}$ et $B_{i,j}=j^t\beta_{i-j}$
pour des suites num\'eriques 
$\alpha_0,\alpha_1,\dots,\beta_0,\beta_1,\dots\in {\mathbb C}^{\mathbb
  N}$, \'etendues \`a des indices 
n\'egatives en posant $\alpha_n=\beta_n=0$ pour tout entier $n<0$.
Dans la suite, on supposera \'egalement $i\geq j$.

La preuve de l'assertion (i) est par r\'ecurrence sur l'exposant
$s$ apparaissant dans les coefficients $A_{i,j}=j^s\alpha_{i-j}$ de 
la matrice $A$. La preuve pour $s=0$ r\'esulte de la d\'efinition
$$(AB)_{i,j}=\sum_k \alpha_{i-k}j^t\beta_{k-j}=j^t\sum_k
\alpha_{i-k}\beta_{k-j}$$
(toutes les sommes sont finies).

Consid\'erons maintenant
$$
(AB)_{i,j}=j^t\sum_k k^s\alpha_{i-k}\beta_{k-j}=C_{i,j}-D_{i,j}$$
avec 
$$C_{i,j}=j^t\sum_k\left(k^s-j^s\right)\alpha_{i-k}\beta_{k-j}$$
et 
$$D_{i,j}=j^{s+t}\sum_k \alpha_{i-k}\beta_{k-j}\ .$$
Les coefficients $D_{i,j}$ sont de la forme voulue. 
La factorisation $(k^s-j^s)=(k-j)\sum_{h=0}^{s-1}k^hj^{s-1-h}$
montre que l'hypoth\`ese de r\'ecurrence s'applique
\`a
$$C_{i,j}=\sum_{h=0}^{s-1}j^{s+t-1-h}\sum_k k^h\alpha_{i-k}\left(
(k-j)\beta_{k-j}\right)$$
ce qui prouve l'assertion (i).

Pour d\'emontrer l'assertion (ii), il suffit de remarquer que le terme
de plus haut
degr\'e $s+t$ dans la preuve de l'assertion (i) correspond \`a
$$\sum_k \alpha_{i-k}\beta_{k-j}=\sum_k \beta_{i-k}\alpha_{k-j}\ .$$
Il se simplifie donc dans le crochet de Lie
$[M(a),M(b)]=M(a)M(b)-M(b)M(a)$.\hfill $\Box$

\begin{rem} \label{autrepreuvelemalg}
(i) On peut aussi prouver le lemme \ref{lemalg} de la mani\`ere
suivante :
Pour $\alpha=\sum_{n=0}^\infty \alpha_nx^n\in{\mathbb C}[[x]]$ 
une s\'erie formelle, notons 
$P_\alpha$ la matrice triangulaire inf\'erieure de coefficients
$(P_\alpha)_{i,j}=\alpha_{i-j}$ (en utilisant 
la convention $\alpha_n=0$ pour
$n<0$) et notons $D$ la matrice diagonale de coefficients diagonaux
$D_{i,i}=i$ pour $i=0,1,2,\dots$.

La matrice $A$ avec coefficients $A_{i,j}=j^s\alpha_{i-j}$ 
consid\'er\'ee dans la preuve du lemme \ref{lemalg} est alors
donn\'ee par $A=P_\alpha D^s$ et le
lemme \ref{lemalg} est une cons\'equence facile du
calcul
$$([P_\alpha,D])_{i,j}=\alpha_{i-j}j-i\alpha_{i-j}=-(i-j)\alpha_{i-j}$$
qui montre l'identit\'e
$$[P_\alpha,D]=-P_{x\alpha'}$$
avec $x\alpha'=\sum_{n=1}^\infty n\alpha_n x^n$ 
pour le crochet de Lie $[P_\alpha,D]=P_\alpha D-DP_\alpha$ de
$P_\alpha$ avec $D$.

\ \ (ii) La partie (i) de la remarque montre que toute matrice
$M\in{\mathcal P}$ est de la forme
$$M=\sum_{n=0}^\infty P_{\alpha(n)}D^n$$
o\`u les s\'eries formelles 
$\alpha(n)=\sum_{n=0}^\infty \alpha(n)_jx^j\in{\mathbb C}[[x]]$ 
v\'erifient la condition suivante : pour tout entier naturel
$k\in{\mathbb N}$ fix\'e,
il existe un entier naturel $N_k\in{\mathbb N}$ tel que
$\alpha(n)_k=0$ pour tout $n> N_k$. 
Il y a alors au plus $N_k$  
contributions non-nulles
\`a $M_{k+j,j}=\sum_{n\in{\mathbb N}}\alpha(n)_{k}j^n$ 
pour $k$ fix\'e et $N_k$ majore le degr\'e du polyn\^ome 
$p_k$ d\'efinie par les coefficients $M_{k+i,i}=p_k(i)$. 
R\'eciproquement, toute 
expression $\sum_{n=0}^\infty P_{\alpha(n)}D^n$ avec
$\alpha(0),\alpha(1),\dots\in{\mathbb C}[[x]]$ des s\'eries formelles 
comme ci-dessus, d\'efinit un
\'el\'ement de ${\mathcal P}$.
\end{rem}

{\bf Preuve du th\'eor\`eme \ref{thmalg}} Le lemme \ref{lemalg} montre
que ${\mathcal P}$ est une alg\`ebre. 
Pour v\'erifier que $\mathcal{SG}$ est une groupe de Lie d'alg\`ebre
$\mathfrak{sg}$, il suffit de 
remarquer que le logarithme matriciel 
$$\hbox{log}(A)=\sum_{n=1}^\infty(-1)^{n+1}\frac{(A-\hbox{id})^n}{n}$$
induit une bijection de $\mathcal{SG}$ sur $\mathfrak{sg}$,
dont l'inverse est l'exponentielle matricielle
$$\hbox{exp}(B)=\sum_{n=0}^\infty\frac{B^n}{n!}\ .$$
Ces propri\'et\'es r\'esultent de l'observation que 
le logarithme de $A\in\mathcal{SG}$ et l'exponentielle de 
$B\in\mathfrak{sg}$ ci-dessus se r\'eduisent \`a des sommes 
finies apr\`es restriction aux sous-matrices finies donn\'ees par 
les $n$ premi\`eres lignes et colonnes. 

L'assertion (ii) est \'evidente pour $\lambda\in\mathbb N$ car
$\tau_\lambda$ efface alors les $\lambda$ premi\`eres lignes et
colonnes d'un \'el\'ement $M\in\mathcal P$. Comme les polyn\^omes sont
des
fonctions analytiques sur $\mathbb C$, le cas g\'en\'eral s'ensuit.
\hfill$\Box$

\begin{rem}
En consid\'erant des matrices index\'ees par ${\mathbb Z}$, on peut
d\'efinir
une alg\`ebre (de Lie) \`a partir de
suites biinfinies de la forme $p=(\dots,0,0,p_N,p_{N+1},\dots)$
avec $N\in{\mathbb Z}$ et $p_i\in{\mathbb C}[u]$, voir la remarque
\ref{remLaurent}.

On peut m\^eme affaiblir cette hypoth\`ese en imposant des conditions
de d\'ecroissance en $n\longrightarrow \pm \infty$ pour les 
polyn\^omes $p_n$.

On peut \'egalement d\'efinir des sous-alg\`ebres (de Lie)
dans ${\mathcal P}$ en imposant des conditions suppl\'ementaires
aux suites de polyn\^omes
$p=(p_0,p_1,\dots)\in{\mathbb C}[u]^{\mathbb N}$
admissibles : on peut par exemple se restreindre aux suites n'ayant
qu'un nombre fini de termes non-nuls. Un autre type de restrictions,
donn\'e par des conditions sur la suite
$\hbox{deg}(p_0),\hbox{deg}(p_1),\dots$ sera consid\'er\'e plus tard.

Toutes les sous-alg\`ebres (de Lie) de ${\mathcal P}$ sont filtr\'ees
sur ${\mathbb N}$. En effet, notons ${\mathcal P}_n\subset {\mathcal
  P}$ les sous-espace vectoriel form\'e de toutes les matrices $M(p)$
avec $p=(p_0=\dots=p_{n-1}=0,p_n,p_{n+1}\dots)\in{\mathbb
  C}[u]^{\mathbb N}$. On a alors ${\mathcal P}_n{\mathcal P}_m\subset
{\mathcal P}_{n+m}$. En particulier,
${\mathcal P}_n$ est un id\'eal bilat\`ere de ${\mathcal P}$. 
\end{rem}

\subsection{Sous-alg\`ebres de degr\'e ${\mathcal B}$ 
dans ${\mathcal P}$}\label{secsousalgB}

Soit ${\mathcal B}\subset{\mathbb N}^{\mathbb N}$ 
un ensemble de fonctions de ${\mathbb N}$ dans
${\mathbb N}$. Nous dirons que ${\mathcal B}$ est {\it additivement
ferm\'e} si pour tout $\alpha,\beta\in {\mathcal B}$
il existe $\gamma\in{\mathcal B}$ satisfaisant
$$\hbox{sup}_{0\leq h\leq k}(\alpha(h)+\beta(k-h))\leq \gamma(k)$$
pour tout $k\in{\mathbb N}$.
On dira qu'une suite de polyn\^omes $p=(p_0(u),p_1(u),\dots)\in{\mathbb
  C}[u]^{\mathbb N}$
est {\it de degr\'e (au plus) ${\mathcal B}$} s'il existe
une fonction $\alpha\in {\mathcal B}$ telle que
$$\hbox{deg}(p_k)\leq \alpha(k)$$ pour tout $k\in {\mathbb N}$.

Soit ${\mathcal B}\subset{\mathbb N}^{\mathbb N}$ 
un ensemble de fonctions ${\mathbb N}
\longrightarrow {\mathbb N}$. Pour $\alpha\in{\mathcal B}$,
notons $\underline{\alpha}$ la fonction d\'efinie par 
$\underline{\alpha}(0)=0$ et $\underline{\alpha}(n)=\alpha(n)$ pour
$n\geq 1$. Notons $\underline{\mathcal B}$ l'ensemble des fonctions
$\underline \alpha$ pour $\alpha\in {\mathcal B}$. Il est 
\'evident que $\underline{\mathcal B}$ 
est additivement ferm\'e si ${\mathcal  B}$ l'est.

On munira les ensembles de fonctions ${\mathbb N}\longrightarrow 
{\mathbb N}$ additivement ferm\'es d'un ordre partiel
en posant ${\mathcal B}_1\leq {\mathcal B}_2$ si toute fonction
$\alpha\in{\mathcal B}_1$ admet une fonction majorante
$\beta\in{\mathcal B}_2$. En d'autres termes, pour tout
$\alpha\in{\mathcal B}_1$, il existe une fonction 
$\beta\in{\mathcal B}_2$ telle que $\alpha(n)\leq \beta(n)$ pour tout 
$n\in{\mathbb N}$.
Par d\'efinition du degr\'e, toute suite $p\in{\mathbb C}[u]^{\mathbb N}$
de degr\'e ${\mathcal B}_1$ est \'egalement de degr\'e
 ${\mathcal B}_2$ si ${\mathcal B}_1,{\mathcal B}_2$ sont
deux ensembles de fonctions additivement ferm\'es tels que 
${\mathcal B}_1\leq {\mathcal B}_2$.

{\bf Exemples} L'ensemble de fonctions 
${\mathcal B}=\{0\}$ r\'eduit \`a la
fonction identiquement nulle est additivement ferm\'e.

L'ensemble ${\mathcal B}$ form\'e de toutes les fonctions constantes
${\mathbb N}\longrightarrow {\mathbb N}$ est additivement ferm\'e.

Pour tout entier $a\geq 0$ et pour tout entier $b\geq 1$,
l'ensemble constitu\'e par la seule fonction $n\longmapsto
a\ n^b$ est additivement ferm\'e.
 
Pour ${\mathcal B}_1,{\mathcal B}_2,\dots$ des ensembles additivement
ferm\'es, l'ensemble des fonctions $\alpha_1+\alpha_2+\dots$
avec $\alpha_i\in {\mathcal B}_i$
est additivement ferm\'es (dans le cas o\`u l'ensemble 
${\mathcal B}_1,{\mathcal B}_2,\dots$ est infini, 
on peut se contenter de sommes sur des sous-ensembles finis).

Pour $p=(p_0,p_1,\dots)\in{\mathbb C}[u]^{\mathbb N}$ 
une suite de polyn\^omes, notons
$p'=(p_0',p_1',\dots)$ la suite des d\'eriv\'ees.
Pour ${\mathcal B}\subset {\mathbb N}^{\mathbb
  N}$ un ensemble additivement ferm\'e, notons ${\mathcal P}^{\mathcal B}$
l'espace vectoriel
contenant toutes les matrices $M(p)$ associ\'ees \`a des suites
de polyn\^omes $p\in{\mathbb C}[u]^{\mathbb N}$
de degr\'e ${\mathcal B}$. Consid\'erons \'egalement
l'espace vectoriel $\mathfrak{sg}^{\mathcal B}\subset
\mathfrak{sg}$ form\'e des matrices
$M(p)$ associ\'ees aux suites 
$p=(p_0,p_1,p_2,\dots)$ avec
$p'=(p_0'=0,p'_1,p_2',\dots)$ de degr\'e ${\mathcal B}$.

\begin{prop} (i) Pour ${\mathcal B}$ un ensemble de fonction ${\mathbb
    N}\longrightarrow {\mathbb N}$ additivement ferm\'e,
l'ensemble des matrices ${\mathcal P}^{\mathcal B}\subset{\mathcal P}$ 
est une sous-alg\`ebre 
et $\mathfrak{sg}^{\mathcal B}\subset\mathfrak{sg}$ est une 
sous-alg\`ebre de Lie.

(ii) Soit ${\mathcal B}\subset {\mathbb N}^{\mathbb N}$ 
additivement ferm\'e. Supposons que toute fonction $\alpha\in{\mathcal
  B}$ admet un majorant $\beta\in{\mathcal B}$ avec $\{\beta\}$ 
additivement ferm\'e. Alors l'exponentielle matricielle 
induit une bijection entre
l'alg\`ebre de Lie $\mathfrak{sg}\cap {\mathcal P}^{\mathcal B}$ et 
le groupe 
des matrices unipotentes $\mathcal{SG}\cap{\mathcal P}^{\mathcal B}$.
\end{prop}

La preuve est une cons\'equence facile du Lemme \ref{lemalg}. 
Elle est laiss\'ee au lecteur.

\begin{rem}
La relation d'ordre sur les sous-ensembles additivement ferm\'es de
${\mathbb N}^{\mathbb N}$ induit
des inclusions au niveau des alg\`ebres (de Lie).

Pour ${\mathcal B}\subset {\mathbb N}^{\mathbb N}$ additivement
ferm\'e, l'ensemble ${\mathcal P}^{\mathcal B}$ est \'egalement une 
alg\`ebre de Lie (pour le crochet usuel $[X,Y]=XY-YX$)
et son alg\`ebre de Lie d\'eriv\'ee $[{\mathcal P}^{\mathcal B},
{\mathcal P}^{\mathcal B}]$ est une sous-alg\`ebre (parfois stricte) de 
$\mathfrak{sg}^{\mathcal B}$.
\end{rem}

{\bf Exemples} L'alg\`ebre ${\mathcal M}^{\{0\}}$ associ\'ee \`a
l'ensemble des fonctions additivement ferm\'e r\'eduit \`a la fonction 
identiquement nulle est l'ensemble 
$$\left(\begin{array}{ccccc}
\alpha_0\\
\alpha_1&\alpha_0\\
\alpha_2&\alpha_1&\alpha_0\\
\alpha_3&\alpha_2&\alpha_1&\alpha_0\\
\vdots  &        &        &        &\ddots\end{array}\right),\alpha_0,
\alpha_1,\dots\in{\mathbb C}$$
des {\it matrices de Toeplitz}
triangulaires inf\'erieures. Une telle matrice est compl\`etement
d\'ecrite par la suite $\alpha_0,\alpha_1,\alpha_2,\dots\in{\mathbb
  C}^{\mathbb N}$ des coefficients de sa premi\`ere colonne
$(\alpha_0,\alpha_1,\dots)^t$ et l'application 
$T(\alpha_0,\alpha_1,\dots)\longmapsto \sum_{n=0}^\infty \alpha_n x^n$ 
qui associe \`a une telle matrice de Toeplitz la s\'erie
g\'en\'eratrice $\sum_{n=0}^\infty \alpha_n x^n$ de sa premi\`ere colonne
est un isomorphisme d'alg\`ebre (pour la structure d'alg\`ebre 
commutative donn\'ee par le produit des s\'eries formelles sur
l'espace vectoriel ${\mathbb C}[[x]]$).

L'alg\`ebre de Lie $\mathfrak{sg}^{\{0\}}$ associ\'ee \`a ${\mathcal
  B}=\{0\}$ est l'espace vectoriel de toutes les matrices triangulaires
inf\'erieures strictes associ\'ees \`a des suites de polyn\^omes
affines. L'alg\`ebre de Lie $\mathfrak{sg}^{\{0\}}$ est donc
l'alg\`ebre de Lie $\mathfrak {si}$ du groupe d'interpolation
sp\'ecial $\mathcal{SI}$ \'etudi\'e dans les chapitres pr\'ec\'edents.

Plus g\'en\'eralement, pour un entier $\lambda\in {\mathbb N}$,
l'ensemble $\{\lambda\hbox{ id}\}$ r\'eduit \`a la fonction lin\'eaire
$n\longmapsto \lambda n$ est additivement ferm\'e. 
On peut donc consid\'erer l'alg\`ebre ${\mathcal P}^{\{\lambda 
\hbox{ id}\}}$. Cette alg\`ebre contient l'id\'eal ${\mathcal S}$
constitu\'e de toutes les matrices de la forme $M(p)$ avec
$p=(p_0=0,p_1,p_2,\dots)\in{\mathbb C}[u]$ o\`u $\hbox{deg}(p_j)<
\lambda j$. La preuve du lemme \ref{lemalg} montre que la 
structure d'alg\`ebre du quotient ${\mathcal M}^{\{\lambda 
\hbox{ id}\}}/{\mathcal S}$ ne d\'epend pas de $\lambda$. L'alg\`ebre
quotient ${\mathcal M}^{\{\lambda 
\hbox{ id}\}}/{\mathcal S}$ s'identifie donc \`a l'alg\`ebre 
${\mathcal M}^{\{0\}}$  (dans laquelle l'id\'eal ${\mathcal S}$ 
se r\'eduit \`a $\{0\}$) des matrices de Toeplitz
consid\'er\'ee auparavant. Remarquons que l'exponentielle matricielle
induit une surjection entre ${\mathcal M}^{\{\lambda 
\hbox{ id}\}}$ et le groupe des \'el\'ements inversibles dans 
${\mathcal M}^{\{\lambda \hbox{ id}\}}$ form\'e de toutes les 
matrices avec diagonale non-nulle dans $\mathcal M^{\{\lambda 
\hbox{id}\}}$. Cette surjection se restreint
en une bijection entre matrices triangulaires inf\'erieures 
strictes dans ${\mathcal M}^{\{\lambda 
\hbox{ id}\}}$ et le sous-groupe ${\mathcal M}^{\{\lambda 
\hbox{ id}\}}\cap \mathcal{SG}$ form\'es des matrices unipotentes
inversibles dans ${\mathcal M}^{\{\lambda 
\hbox{ id}\}}$. La bijection inverse est \'evidemment donn\'ee par le
logarithme matricielle (qui converge toujours coefficient par 
coefficient car $G\in\mathcal{SG}$ est une matrice 
triangulaire inf\'erieure unipotente).

Le groupe sp\'ecial d'interpolation
$\mathcal{SI}=\{\hbox{exp}(A)\ \vert\ A\in\mathfrak{si}\}$, est
un sous-groupe de l'alg\`ebre ${\mathcal M}^{
\{\hbox{id}\}}$ associ\'ee \`a $\lambda=1$. Ajoutons que 
l'alg\`ebre  ${\mathcal M}^{\{2
\hbox{ id}\}}$ (qui contient \'evidemment ${\mathcal M}^{\{
\hbox{id}\}}$ comme sous-alg\`ebre) contient \'egalement des 
\'el\'ements int\'eressants comme par exemple une matrice triangulaire
inf\'erieure dont les coefficients sont donn\'es par les
nombres de Stirling  (de premi\`ere ou deuxi\`eme esp\`ece).


Institut Fourier,
Laboratoire de Math\'ematiques,
UMR 5582 (UJF-CNRS),
100, rue des Math\'ematiques,
BP 74,
38402 St MARTIN D'H\`ERES Cedex, France

Adresse courriel: Roland.Bacher@ujf-grenoble.fr


\begin{thebibliography}{99}

\bibitem{CN} N.T.~Cameron, A.~Nkwanta,
{\it On Some (Pseudo) Involutions in the
    Riordan Group}, J. of Int. Seq., Vol 8 (2005), Article 05.3.7.

\bibitem{EIS} N.J.A.~Sloane, {\it The On-Line 
Encyclopedia of Integer
Sequences}, http://www.research.att.com/~njas/sequences/index.html

\bibitem{KPP} A.G.~Kuznetsov, I.M.~Pak, A.E.~Postnikov, {\it
Increasing trees and alternating permutations},	
Russ. Math. Surv. {\bf 49}, No.6, 79-114 (1994); 
traduit de : Usp. Mat. Nauk {\bf 49}, No.6, 79-110 (1994).

\bibitem{Labelle} G.~Labelle, {\it Sur l'inversion et l'iteration continue
des series formelles}, Europ. J. Combin., {\bf 1} (1980), 113--138.

\bibitem{MRSV} D.~Merlini, D.G.~Rogers, R.~Sprignoli, M.C.~Verri, {\it
    On some alternative characterizations of Riordan arrays},
  Can. J.Math. {\bf 49}, No. 2, 301-320 (1997). 


\bibitem{R} D.~G.~ Rogers, {\it Pascal triangles, Catalan numbers and
    renewal arrays}, Discr. Maths. {\bf 22} (1978), 301-310.

\bibitem{SGWW} L.~W.~Shapiro, S.~Getu, W.~J.~Woan, L.~C.~Woodson, {\it
The Riordan group}, Discrete Appl. Math.  {\bf 34}  (1991), 
no. 1-3, 229-239.

\bibitem{Spr} R.~Sprugnoli, An Introduction to Mathematical Methods in
  Combinatorics, preprint.

\bibitem{St1} R.P.~Stanley, Enumerative Combinatorics, Vol. 1, 
Cambridge University Press, 1997.



\end{thebibliography}
\end{document}